\documentclass[a4paper,10pt]{article} \usepackage{amsmath,amssymb,amsthm}
\usepackage[english]{babel}

\newtheorem{theo}{Theorem}[section] \newtheorem{defi}[theo]{Definition}
\newtheorem{lemm}[theo]{Lemma} \newtheorem{prop}[theo]{Proposition}
\newtheorem{coro}[theo]{Corollary}

\newcommand{\Na}{\mathbb N}                   % N des entiers

\newcommand{\Ra}{\mathbb R}                   % R des reels

\newcommand{\Ca}{\mathbb C}                   % C des complexes

\newcommand{\scal}[1]{\langle #1 \rangle}

\newcommand{\finpreuve}{\hfill $\Box$}

\newcommand{\name}{$\underline{\qquad \qquad}$}

\newcommand{\refe}[1]{\ref{#1}} \newcommand{\reff}[1]{(\ref{#1})}

\begin{document}

\author{  Jean-Marc
Bouclet\footnote{Jean-Marc.Bouclet@math.univ-lille1.fr}\\
 \\ Universit\'e  de Lille  1 \\ Laboratoire Paul Painlev\'e \\ UMR  CNRS
8524,  \\ 59655 Villeneuve  d'Ascq }
\title{{\bf \sc Semi-classical functional calculus on 
    manifolds with ends and weighted  $L^p$ estimates}}

\maketitle

\begin{abstract} For a class of non compact Riemannian manifolds with ends $ ({\mathcal M},g) $, we give pseudo-differential expansions of bounded functions  of the semi-classical Laplacian $ h^2 \Delta_g $, $ h \in (0,1] $. We then study related $ L^p $ boundedness properties and show in particular that, although
 $ \varphi (-h^2 \Delta_g) $ is not bounded on $ L^p ({\mathcal M},dg) $ in general, it is always bounded on suitable weighted $ L^p $ spaces. 
 %in terms of properly supported $h$-pseudo-differential operators
\end{abstract}

$$  $$

\section{Introduction and Results}

In this paper we  describe semi-classical expansions of functions of the
Laplacian on a class of non compact manifolds of bounded geometry. We also
derive certain (weighted) $ L^p \rightarrow L^p $  boundedness properties of
such operators. Further applications to Littlewood-Paley decompositions
\cite{BoLP} and Strichartz estimates \cite{BoSH} will be published
separately. Needless to say,  the range of applications of the present functional
calculus goes beyond Strichartz estimates; there are many problems which
naturally involve  spectral cutoffs at high frequencies in linear and non linear
PDEs (Littlewood-Paley decompositions, paraproducts)  or in spectral theory (trace formulas).

Consider a  non compact  Riemannian manifold $ ({\mathcal M},g) $ with ends, ie whose model at infinity is a product $ (R,+\infty) \times S $ with metric  $ g = d r^2 + d \theta^2 / w(r)^2 $, where $ R \gg 1 $, $(S,d\theta^2)$  a compact Riemannian manifold
and $w(r)$ a bounded positive function. For instance,  $ w (r) = r^{-1} $ corresponds to conical ends, $ w (r) = 1 $ to cylindrical ends and $ w(r) = e^{-r} $ to hyperbolic ends. We actually consider more general metrics (see Definition \ref{defmanifold} below for precise statements)
but these are the typical examples we have in mind.
If $ \Delta_g $ denotes the Laplacian on $ {\mathcal M} $ and $ \varphi $ is a symbol of negative order, we are interested in decompositions of the form
\begin{eqnarray}
 \varphi (-h^2 \Delta_g) = {\mathcal Q}_N (\varphi,h) + h^{N+1} {\mathcal R}_N (\varphi,h), \qquad h \in (0,1 ] , \label{ansatzvague}
\end{eqnarray}
where $ N \geq 0$ is fixed and arbitrary, $ {\mathcal Q}_N (\varphi,h) $ has an expansion in powers of $h$ in terms $h$-pseudo-differential operators and $ h^{N+1} {\mathcal R}_N (\varphi,h) $ is a 'nice' remainder.  We  recall that, for such semi-classical expansions, even the case of $ \varphi \in C_0^{\infty}(\Ra) $ is of interest, by opposition to the classical case ($h = 1$) where $ C_0^{\infty} $ functions of $ \Delta_g $ are often treated as negligible operators.

There is a large literature devoted to 
the pseudo-differential analysis of functions of closed operators on manifolds so we only give  references
which are either classical or close to our framework. For $ h =1 $, the case
of compact manifolds 
(ie, essentially, the local interior case) was considered by Seeley \cite{Seel1} (see
also  \cite[pp. 917-920]{Seel2}). For boundary
value problems, we refer to \cite{Seel2,Grub1} and
for non compact or singular manifolds to \cite{Schul1,ALNV1}.
%\cite[?]{Schul2}. 
We also quote \cite{CGT1,Tayl0,Kord1} where  general manifolds of bounded
geometry are studied in connection with the problem of the $ L^p \rightarrow L^p  $
boundedness of functions the Laplacian (to which we come back below). The semi-classical case is treated for
very general operators on $ \Ra^n $  in \cite{HeRo1,Robe1,DiSj1} and 
in  \cite{BGT1} for a compact manifold. Besides, one of our initial
motivations is to extend  the functional calculus used
in \cite{BGT1} to non compact manifolds and thus to provide a convenient 
tool to prove Strichartz estimates,  as for instance in \cite{HTW0,BoTz1}.

Although the general picture is quite clear, at least from the $ L^2 $ point
of view, the problem of getting expansions of the 
form (\ref{ansatzvague})
% of pseudo-differential functional calculus on non compact manifolds 
requires some care. By opposition to the compact case (or to $ \Ra^n $ for
uniformly elliptic operators), one
has to take into account certain off diagonal effects possibly leading to the
unboundedness of the operators on $ L^p ({\mathcal M},dg) $, when $ p \ne 2 $,
if $dg$ denotes the Riemannian measure. 

By considering properly supported
operators, namely with kernels supported close to the diagonal of $ {\mathcal
  M} \times {\mathcal M} $, we may insure that the principal part of the
expansion $ {\mathcal Q}_N (\varphi,h) $ is bounded on $ L^p ({\mathcal M},dg)
$, for all $ p \in [1,\infty] $, uniformly with respect to $h$. However, the
boundedness of the remainder $ {\mathcal R}_N (\varphi,h) $ on $ L^p
({\mathcal M},dg) $ remains  equivalent to the one of the full operator $
\varphi (-h^2 \Delta_g) $ and it is well known that the latter may fail for
non holomorphic $ \varphi $, as
first noticed by Clerc and Stein \cite{ClSt1} for symmetric spaces.
The latter question is treated (with $h=1$)   for a large class
of manifolds by Taylor  in \cite{Tayl0} (see also the
references therein and the extension \cite{Kord1} to systems of properly
supported operators). Taylor proves  that, if $A$ denotes the bottom of the spectrum of $ - \Delta_g $ and $ L = (- \Delta_g -A)^{1/2} $, the boundedness of $ \varphi (L) $ on $ L^p ({\mathcal M},dg) $ is guaranteed if 
 $ \varphi $ is even and holomorphic in a strip of width at least $\kappa|1/p-1/2|$, 
with $ \kappa $ the exponential rate of the volume growth of
balls. This is typically relevant in the hyperbolic case. To illustrate this
fact (as well as some of our results), we recall  a  short proof of the $ L^p
$-unboundedness  of $ (z-\Delta_{{\mathbb H}^n})^{-1} $ in  Appendix
\ref{contreexemple}, 
$ \Delta_{{\mathbb H}^n} $ being the Laplacian on the hyperbolic space.
 
%Such non $ L^p \rightarrow L^p  $ boundedness properties are closely related
%to the fact that pseudo-differential operators do not preserve exponential
%decays in general. Exponential decay may however be preserved but under certain analyticity conditions on the
%symbols. See for instance \cite{Schot} for results in this direction.

In summary, our first goal is to provide a fairly explicit and precise description of expansions of the form (\ref{ansatzvague}). For $h=1$, this result is essentially contained in \cite{CGT1,Tayl0} but we feel that it is worth giving complete proofs for the semi-classical case too, first because we shall use it extensively in subsequent papers and second because of the subtleties due to $ L^p $-unboundedness.

Our second point is to prove weighted $L^p $ estimates on $ {\mathcal R}_N
(\varphi,h) $ or, equivalently, on the resolvent $ (z-\Delta_g)^{-1}
$. The basic strategy is to use the expansion (\ref{ansatzvague}) to get $ L^2 $ estimates on commutators of the resolvent with natural first order
differential operators and show that $
(z-\Delta_g)^{-1}  $ is a pseudo-differential operator, using the Beals criterion. At this stage, the
meaning of pseudo-differential operator is rather vague but we emphasize
that the point  is not (only) to control the singularity of the kernel
close to the diagonal but also the decay far from the diagonal. As a
consequence of this analysis,
 we  obtain in particular that, although $ (z-\Delta_g)^{-1} $ is not necessarily bounded on $ L^p ({\mathcal M},dg) $, we always have
$$ || w(r)^{\frac{n-1}{p} - \frac{n-1}{2}} (z - \Delta_g)^{-1} w(r)^{\frac{n-1}{2} - \frac{n-1}{p}} ||_{L^p ({\mathcal M},dg) \rightarrow L^p ({\mathcal M},dg)} < \infty, $$
for all $ p \in (1,\infty) $ and $ z \notin \mbox{spec}(\Delta_g) $. More
generally, if $ W $ is a temperate weight (see Definition
\ref{definitionpoidstempere} below), we have
$$  || W(r)^{-1} w(r)^{\frac{n-1}{p} - \frac{n-1}{2}} (z - \Delta_g)^{-1}
w(r)^{\frac{n-1}{2} - \frac{n-1}{p}} W(r) ||_{L^p ({\mathcal M},dg)
  \rightarrow L^p ({\mathcal M},dg)} < \infty .  $$
  This
works in particular for the hyperbolic case where $   (z - \Delta_g)^{-1} $ is
not bounded on $ L^p ({\mathcal M},dg) $ in general.  In the conical case, 
or more generally if $ w $ itself is a temperate weight,  we   recover the natural (unweighted)
boundedness on $ L^p ({\mathcal M},dg) $ by choosing $ W = w^{\frac{n-1}{p} -
  \frac{n-1}{2}}   $. The latter boundedness can be seen as a consequence of
\cite{Tayl0} since, if $w$ is temperate, the volume growth of balls is
polynomial. The above estimates are therefore complementary to the results of
\cite{Tayl0}: if $z$ is too close to the spectrum of the Laplacian, $
(z-\Delta_g)^{-1}  $ is maybe not bounded on $ L^p =  L^p ({\mathcal M},dg)  $ but it
is bounded if we accept to replace $ L^p $ by weighted $L^p$ spaces. Furthermore,
these weighted spaces  are natural since they contain $ L^p $ itself when $w$ is
temperate (ie essentially if $w^{-1}$ is of polynomial growth).

%However, we emphasize that we obtain a stronger result: we show that the
%resolvent belongs to an explicit class of 
%pseudo-differential operators from which weighted  $ L^p $ bounds are rather easily deduced.

\bigskip

Let us now state our results precisely.

 In the sequel $ {\mathcal M} $ will be a smooth manifold
of dimension $ n \geq 2$, without boundary and satisfying the following definition.
\begin{defi} \label{defmanifold} The manifold $ ({\mathcal M},g ) $ is called almost asymptotic
if  there exists a compact set $ {\mathcal K} \Subset {\mathcal M}
$, a real number $ R  $, a compact manifold $ S $, a function $ r \in C^{\infty}({\mathcal M},\Ra) $ and a function $
w \in C^{\infty}( \Ra , (0,+\infty) ) $  with the following
properties: 
\begin{enumerate}
\item{ $r$ is a coordinate near $ \overline{{\mathcal M} \setminus {\mathcal K}}  $ such that
$$ r (\underline{x}) \rightarrow + \infty, \qquad \underline{x} \rightarrow \infty , $$ }
\item{ there is a diffeomorphism
\begin{eqnarray}
  {\mathcal M} \setminus {\mathcal K}  \rightarrow (R,+\infty) \times S, \label{diffeomorphismeabstrait}
\end{eqnarray}
through which  the metric reads, in local coordinates,
\begin{eqnarray}
 g & = & G_{\rm unif} \left(
r , \theta , d r ,  w(r)^{-1} d\theta \right)
\label{pourlelaplacien} 
%& = & G_{11}(r,\theta)dr^2 + 2 G_{1 k+1}(r,\theta)dr
%w(r)^{-1}d\theta_k + G_{j+1 k+1}(r,\theta)w(r)^{-2} d \theta_j d
%\theta_k
\end{eqnarray}
  with 
$$ G_{\rm unif}(r,\theta,V) := \sum_{1 \leq j,k \leq n} G_{jk}(r,\theta) V_j V_k , \qquad V = (V_1, \ldots,V_n)\in \Ra^n , %+ 2 G_{1
%k+1}(r,\theta) V \eta_{k}  + G_{j+1 k+1}(r,\theta) \eta_j \eta_k
$$
if  $ \theta = (\theta_{1}, \ldots , \theta_{n-1}) $ are local
coordinates on $S$.}
\item{ The symmetric matrix $ (G_{jk}(r,\theta))_{1 \leq j,k \leq n} $ has smooth coefficients
such that, locally uniformly with respect to $ \theta $,
% for $ \rho \in \Ra $ and $ \eta = (\eta_1,\ldots,\eta_{n-1}) \in \Ra^{n-1}
%$. We assume that  the functions $G_{jk} = G_{kj}$ are smooth and
%satisfy\footnote{ Throughout the paper $ a (x) \lesssim b(x), \ x
%\in E $, will mean that there exists $ C
%> 0 $ such that, for all $ x \in E  $, $ a (x) \leq C b
%(x) $. The notation $ a (x) \thickapprox b (x), \ x \in E $, will
%mean that $ a (x) \lesssim b(x)  \ \mbox{and} \  b (x) \lesssim a
%(x) , \ x \in E $.}, 
\begin{eqnarray}
 \left| \partial_r^j \partial_{\theta}^{\alpha} G_{jk}(r,\theta) \right| & \lesssim & 1, \qquad \qquad r
> R,  \label{uniformegunif} 
\end{eqnarray}
and is uniformly positive definite in the sense that, locally uniformly in $ \theta $,
\begin{eqnarray}
G_{\rm unif}(r,\theta,V) & \thickapprox & |V|^2 , \qquad r
> R, \ V \in \Ra^{n}. \label{bornespectre}
\end{eqnarray}}
\item{ The function $ w $ is smooth and satisfies,
for all $ k \in \Na $,
\begin{eqnarray}
 w (r) & \lesssim & 1 ,  \label{sanscusp} \\
w (r) / w (r^{\prime}) & \thickapprox & 1 , \qquad \mbox{if} \ \ |r-r^{\prime}| \leq 1 \label{diag}  \\
\left|  d^k w (r) / dr^k \right| & \lesssim & w (r),
\label{cinfiniborne}
\end{eqnarray}
for  $ r,r^{\prime} \in \Ra $. }
\end{enumerate}
\end{defi}

Note that (\ref{diag}) is equivalent to the fact that, for some
$ C > 0 $,
$$ C^{-1} e^{-C|r-r^{\prime}|} \leq \frac{w (r)}{  w (r^{\prime}) } \leq C e^{C|r-r^{\prime}|} . $$ 
In particular, this implies that $ w
(r) \gtrsim e^{-C |r|} $.

Asymptotically conical
manifolds, for which $ g = d r^2 + r^2 g_S (r,\theta,d\theta) $ (near infinity), or
  asymptotically hyperbolic manifolds  for which
$ g = d r^2 + e^{2 r} g_S (r,\theta,d\theta) $,
with $ g_S (r, \theta,d \theta) $  a metric on $ S $ depending
smoothly on $r$, satisfy our definition. More precisely, for such asymptotic structures
one usually requires that $ g_S (r,\theta,d \theta) $ is a small
perturbation of a metric $ g_S^{\infty}(\theta,d\theta) $ in the
sense that $   g_S (r,\theta,d\theta) -
g_S^{\infty}(\theta,d\theta)  \rightarrow 0 $ as $ r \rightarrow
\infty $. See for instance \cite{Melr0} for more precise
statements. Here we do not require such a condition which is the
reason why we use the terminology {\it almost} asymptotic.

\medskip

{\it Atlas and partition of unit.} We now specify an
atlas on $ {\mathcal M} $. The diffeomorphism  (\ref{diffeomorphismeabstrait}) is of the form $
\Psi : {\mathcal M} \setminus {\mathcal K} \rightarrow (R, +
\infty) \times { S } $ with
$$ \Psi (\underline{x}) = (r(\underline{x}),\pi_{ S}(\underline{x})), \qquad  \underline{x}
 \in {\mathcal M} \setminus
{\mathcal K} ,  $$ where $ \pi_{ S} : {\mathcal M} \setminus
{\mathcal K} \rightarrow { S} $ is the "projection on the manifold
at infinity" and $ r   $ the "radial coordinate" used in
Definition \ref{defmanifold}. Thus, if we consider a
chart on $S$, 
$$ \psi_{\iota} : U_{\iota} \subset { S}
\rightarrow V_{\iota} \subset \Ra^{n-1} , $$ 
with $ \psi_{\iota} (\underline{y}) = (\theta_1(\underline{y}),\ldots,\theta_{n-1}(\underline{y})) $, then the open sets
 \begin{eqnarray}
     {\mathcal U}_{\iota} =  \Psi^{-1} \left( (R,+\infty) \times U_{\iota} \right) \subset {\mathcal M}  ,
  \qquad {\mathcal V}_{\iota} =  (R,+\infty) \times V_{\iota} \subset \Ra^n ,
  \label{ouverts}
  \end{eqnarray}
and the map
\begin{eqnarray*}
 \Psi_{\iota} : {\mathcal U}_{\iota}  
\rightarrow  {\mathcal V}_{\iota} , \qquad \mbox{with} \ \
\ \Psi_{\iota}(\underline{x}) & = & (r(\underline{x}), \psi_{\iota}
\circ \pi_S (\underline{x})) \\
& = & \left( r(\underline{x}), \theta_1 ( \pi_S (\underline{x})) , \ldots, \theta_{n-1}
 (\pi_S (\underline{x})) \right) ,
\end{eqnarray*}
 define a coordinate chart on $ {\mathcal M} \setminus
{\mathcal K}  $. With a standard abuse of notation, we will denote for simplicity these coordinates $ (r,\theta_1,\ldots,\theta_{n-1}) $
or even $ (r,\theta) $.
\begin{defi} We call $ {\mathcal U}_{\iota} $ a coordinate
  patch at infinity and  the triple $ ({\mathcal U}_{\iota},{\mathcal
V}_{\iota},\Psi_{\iota}) $ a  chart at infinity.
\end{defi}
Since $ { S}  $ is compact, there is a finite set $
I_{\infty} $  such that the family $
({\mathcal U}_{\iota}, {\mathcal V}_{\iota},
\Psi_{\iota})_{\iota \in I_{\infty}} $ is an atlas on $
{\mathcal M} \setminus {\mathcal K} $. Choosing another finite
collection of coordinate charts for a neighborhood of $ {\mathcal
K} $, which we denote\footnote{we keep the  notation $ {\mathcal
U}_{\iota}, {\mathcal V}_{\iota}, \Psi_{\iota} $
  but, of course, the corresponding new $ {\mathcal U}_{\iota} $ and
 $ {\mathcal V}_{\iota} $ are not defined by  (\ref{ouverts}). In the core of the paper, there should be anyway no confusion for
 we shall work almost only on $ {\mathcal M} \setminus {\mathcal K} $.} by $ ({\mathcal U}_{\iota}, {\mathcal
V}_{\iota}, \Psi_{\iota})_{\iota \in I_{\rm comp}} $ for
some finite set $ I _{\rm comp} $,  we get a finite atlas on
$ {\mathcal M} $ by considering $ ({\mathcal U}_{\iota},
{\mathcal V}_{\iota}, \Psi_{\iota})_{\iota \in I} $ with
$$ I = I_{\infty} \cup I_{\rm comp} . $$ In
particular, we can find a finite partition of unit
\begin{eqnarray}
 \sum_{\iota \in I} f_{\iota} = 1 \qquad \mbox{on} \ {\mathcal M} ,  \label{partition}
\end{eqnarray}
such that, for all $ \iota \in I $, $ f_{\iota}$ is
supported in $ {\mathcal U}_{\iota} $. We also set
\begin{eqnarray}
 \chi_{\iota} = f_{\iota} \circ \Psi_{\iota}^{-1} . \label{notationcutoff1}
\end{eqnarray}
 If $ {\mathcal U}_{\iota}
$ is a patch at infinity, we can assume that $ f_{\iota} $ is such that
\begin{eqnarray}
  \chi_{\iota}(r,\theta)  = \varrho (r) \kappa_{\iota}(\theta) ,
 \label{notationcutoff2}
\end{eqnarray}
for some smooth functions $ \varrho $ and $ \kappa_{\iota} $ such
that, for some $ R^{\prime} > R $,
\begin{eqnarray}
 \varrho (r) = 1 \ \ \mbox{for} \   r \gg  1, \ \  \mbox{supp} \
\varrho \subset [R^{\prime},+\infty) , \qquad  \kappa_{\iota} \in
C_0^{\infty} ( V_{\iota} ). \label{choixtroncature}
\end{eqnarray}

\medskip

{\it Differential operators on $ {\mathcal M} $.} We first compute
the Laplacian $ \Delta_g $ in a chart at infinity. Let us define $
\partial_1^{w}, \ldots,
\partial_n^{w} $ by
$$ \partial_1^w = \partial_r , \qquad \partial_2^w = w(r) \partial_{\theta_1},
\ \ldots \ , \ \partial_n^w = w(r) \partial_{\theta_{n-1}}^w . $$
We also  set  $ (G^{jk})_{1 \leq j,k \leq n} := (G_{jk})^{-1}_{1
\leq j,k \leq n} $ and $ \mbox{det} \ G_{\rm unif} := \mbox{det}
(G_{jk}) $ (see (\ref{pourlelaplacien})).  
We then have
\begin{eqnarray}
\Delta_g = (\mbox{det} \ G_{\rm unif})^{-1/2} \partial_j^w G^{jk} (
\mbox{det} \ G_{\rm unif} )^{1/2} \partial_k^w + (1-n)
\frac{w^{\prime}(r)}{w(r)} G^{1k} \partial_k^w ,
\label{laplacienexplicite}
\end{eqnarray}
using the summation convention for $ j , k \geq 1 $.
 This formula motivates the introduction of the following class of differential operators.

\begin{defi} For $ m \in \Na $,  $ \emph{Diff}_w^m ({\mathcal M}) $
is the space of differential operators $ P $ of order $ \leq m $,
acting on functions on $ {\mathcal M} $, such that, for any chart
at infinity $ ({\mathcal U}_{\iota}, {\mathcal V}_{\iota},
\Psi_{\iota}) $,
\begin{eqnarray}
  \Psi_{\iota *} P \Psi_{\iota}^* = \sum_{ k + |\alpha| \leq m } a_{k \alpha}^{\iota}(r,\theta) \left(
   w(r)D_{\theta} \right)^{\alpha} D_r^k  , \label{notationlocalezero}
\end{eqnarray}
with
$$ \partial_r^j
\partial_{\theta}^{\beta} a_{k \alpha}^{\iota} \in L^{\infty}\left(
(R,+\infty) \times K_{\iota} \right)  , $$ for all $ j , \beta $ and all $
K_{\iota} \Subset V_{\iota} $. Here we used
the standard notation $ \Psi_{\iota}^* u = u \circ \Psi_{\iota}
$ and $ \Psi_{\iota * } v = v \circ \Psi_{\iota}^{-1} $.
\end{defi}
By  (\ref{uniformegunif}), (\ref{cinfiniborne})  and  
(\ref{laplacienexplicite}),
 we  see that $ - \Delta_g \in \mbox{Diff}_w^2
({\mathcal M}) $ and that its principal symbol takes the following
form in $ {\mathcal V}_{\iota} $, for $ \iota \in
I_{\infty} $,
\begin{eqnarray}
 p_2^{\iota} (r,\theta,\rho,w(r)\eta) = G^{11}(r,\theta) \rho^2 +
 2 G^{1k}(r,\theta)\rho w(r)\eta_k + G^{j k} (r,\theta) w(r)^2
 \eta_j\eta_k ,
\label{symbolprincipalinfini}
\end{eqnarray}
using the summation convention for $ j,k \geq 2 $. Here and below $ \rho $ and $ \eta $ denote respectively the dual variables to $
r $ and $ \theta $. If $ \iota \in I_{\rm comp}  $, the principal
symbol of $ - \Delta_g $ in $ {\mathcal V}_{\iota} $ takes the
standard form
\begin{eqnarray}
 p_2^{\iota}(x,\xi) = g^{jk}(x) \xi_j \xi_k \label{symbolprincipalcompact}
\end{eqnarray}
  for some smooth $
(g^{jk}(x)) $ such that $ g^{jk}(x) \xi_j \xi_k \gtrsim |\xi|^2 $ for $ \xi \in \Ra^n $ locally uniformly with respect to $x$.

\medskip

\noindent {\bf Remark.} Recall that, if $ \iota \in I_{\infty} $, the principal of $ - \Delta_g $ is given by (\ref{symbolprincipalinfini}) but not by $ p_2^{\iota} $ itself (see the factor
$ w (r) $ in the left hand side of (\ref{symbolprincipalinfini})). This  notation (which is perhaps confusing) will be convenient to state Theorem \ref{theocalcfon1}.

\bigskip

 {\it Lebesgue spaces.} We now describe volume densities. In
coordinates $ (r , \theta) $ at infinity, the Riemannian volume density associated
to $g$,  denoted by $dg $, reads
\begin{eqnarray}
 d g = w (r)^{1-n} (\mbox{det} \ G_{\rm unif}(r,\theta))^{1/2} dr
d\theta, \label{mesure}
\end{eqnarray}
 where, for $ r > R $ and
  locally uniformly with respect to $
\theta $,  (\ref{bornespectre}) yields
\begin{eqnarray}
 \mbox{det} \ G_{\rm unif}(r,\theta) \approx 1 . \label{pourequivalenceLpMRn}
 \end{eqnarray}
Define another
density $ \widetilde{dg} $ on $ {\mathcal M} $ by
\begin{eqnarray}
 \widetilde{dg} = w^{n-1}(r) d g ,  \label{mesuretilde}
\end{eqnarray}
we then have
\begin{eqnarray}
 L^p ({\mathcal M},dg) = w^{\frac{n-1}{p}}(r) L^p ( {\mathcal M},
\widetilde{dg} ) , \qquad p \in [1,\infty). \label{weightlp}
\end{eqnarray}
The map $ u \mapsto w(r)^{(n-1)/2} u $ is unitary from
$ L^2 ({\mathcal M},\widetilde{dg}) $ on $ L^2 ({\mathcal M},dg) $
and the operator
\begin{eqnarray}
  \widetilde{\Delta}_g :=  w(r)^{\frac{1-n}{2}} \Delta_g
w(r)^{\frac{n-1}{2}} , \label{laplacientilde}
\end{eqnarray}
 is symmetric
on $ C_0^{\infty} ({\mathcal M}) $ with respect to $ \widetilde{ d
g } $.\footnote{$ \Delta_g $ and $ \widetilde{\Delta}_g $ are
respectively essentially self-adjoint on $ L^2 ({\mathcal M},dg) $
and $ L^2 ({\mathcal M},\widetilde{dg}) $ from $
C_0^{\infty}({\mathcal M}) $ and thus unitarily equivalent.} By (\ref{cinfiniborne}), we  have 
$$ \widetilde{\Delta}_g \in \mbox{Diff}_w^2({\mathcal M}) . $$

Note that,  for all $ \iota \in I_{\infty} $ and all $ K_{\iota} \Subset V_{\iota} $ (see (\ref{ouverts})), we have the equivalence of norms
\begin{eqnarray}
|| u ||_{L^p ({\mathcal M},\widetilde{dg})} \approx || u \circ \Psi_{\iota}^{-1} ||_{L^p (\Ra^n,drd\theta)}, \qquad \mbox{supp}(u \circ \Psi_{\iota}^{-1}) \subset (R,+\infty) \times K_{\iota} , \label{equivalenceLpMRn}
\end{eqnarray}
for $ p \in [1,\infty] $. This is a simple consequence of  (\ref{pourequivalenceLpMRn}). On compact subsets, the same equivalence holds trivially. For the measure $ dg $, we have, if $ \iota \in I_{\infty} $,
\begin{eqnarray}
\left| \left| u   \right| \right|_{L^p
({\mathcal M}, dg)}  \approx   \left| \left|
w^{(1-n) /p}(r) u \circ \Psi_{\iota}^{-1} \right| \right|_{L^p (\Ra^n ,drd\theta)}   ,
\qquad \mbox{supp}(u \circ \Psi_{\iota}^{-1}) \subset (R,+\infty) \times K_{\iota} . \label{equivalenceLpapoids}
\end{eqnarray}

\bigskip

{\it Pseudo-differential operators.} We now define a
 class of semi-classical pseudo-differential
operators associated to the partition of unit (\ref{partition}). We will choose symbols
$$ a^{\iota} \in S^{m}_{\iota} ({\mathcal V}_{\iota} \times \Ra^n) , $$
where $ {\mathcal V}_{\iota} \subset \Ra^n $ is  defined by (\ref{ouverts}) if $ \iota \in I_{\infty} $. By definition, this means, if $ \iota \in I_{\infty} $, that for all $ K_{\iota} \Subset V_{\iota} $,
$$  | \partial_r^j \partial_{\theta}^{\alpha} \partial_{\rho}^k \partial_{\eta}^{\beta} a^{\iota} (r,\theta,\rho,\eta) | \leq C (1+|\rho|+|\eta|)^{m-k-|\beta|} , \qquad r > R ,\ \ \theta \in K_{\iota}, \ \rho \in \Ra, \ \eta \in \Ra^{n-1} , $$
and, if $ \iota \in I_{\rm comp} $, that  for all $ {\mathcal K}_{\iota} \Subset {\mathcal V}_{\iota} $,
$$ |  \partial_{x}^{\alpha} \partial_{\xi}^{\beta} a^{\iota} (x,\xi) | \leq C (1+|\xi|)^{m-|\beta|} , \qquad x \in {\mathcal K}_{\iota}, \ \xi \in \Ra^{n}. $$ 
%The latter corresponds to the usual class $ S^m_{\rm loc} ({\mathcal V}_{\iota} \times \Ra^n) $. The difference with $ \iota \in I_{\infty} $ is that %we   require the uniformity of the estimates with respect to $ r$. 
In both cases, the topology of $ S^m_{\iota}({\mathcal V}_{\iota} \times \Ra^n) $
is given by the best constants $ C $ which define semi-norms.

We basically would like to use operators of the form
$$ a^{\iota} (r,\theta,h D_r, h w(r)
D_{\theta}) \chi_{\iota} , \qquad \mbox{if} \ \iota \in I_{\infty}, $$
(see (\ref{pseudo}) below) and
$$ a^{\iota} (x, h D_x) \chi_{\iota} , \qquad \mbox{if} \ \iota \in I_{\rm comp} , $$
where $ \chi_{\iota} $ is defined by (\ref{notationcutoff1}) and $h \in (0,1 ] $ is the semi-classical parameter. Actually,
we need to consider properly supported operators
so we construct first suitable cutoffs near the diagonal. Choose  a function $ \zeta \in
C_{0}^{\infty}(\Ra^{n}) $ and $ \varepsilon > 0 $  such that
\begin{eqnarray} \zeta (x) = 1  \ \mbox{ for } \  |x| \leq \varepsilon
, \qquad \zeta (x) = 0  \ \mbox{ for } \ |x| > 2 \varepsilon .
\label{zeta}
\end{eqnarray}
For $ \iota \in I_{\infty} $, the  function 
\begin{eqnarray}
  \chi^{\zeta}_{\iota}(r,\theta,r^{\prime},\theta^{\prime})
:= \chi_{\iota}(r^{\prime},\theta^{\prime}) \zeta \left(
  (r,\theta)-(r^{\prime},\theta^{\prime}) \right), \label{cutoffdiagonal}
\end{eqnarray}
is smooth on $ \Ra^{2n} $  and, 
%we may assume that it is
%supported in $ {\mathcal  V}_{\iota} \times {\mathcal V}_{\iota}
%$ (see (\ref{ouverts}). More precisely, 
if $ K_{\iota} \Subset V_{\iota} $ is an arbitrarily small neighborhood of $ \mbox{supp}(\kappa_{\iota}) $ (see \ref{choixtroncature})), we may choose $ \varepsilon $ small enough such that
\begin{eqnarray}
\mbox{supp}(\chi^{\zeta}_{\iota}) \subset  \left( (R,+\infty) \times K_{\iota} \right)^2 . \label{localisationnoyauoperateur}
\end{eqnarray}
 Proceeding similarly 
 %for those $ f_{\iota} $ which are
%compactly supported, ie 
for $ \iota \in I_{\rm comp} $, we obtain a family of functions $ (\chi_{\iota}^{\zeta})_{\iota \in I} $ supported close to the
diagonal of $ \Ra^{2n} $, with also $ \mbox{supp}( \chi_{\iota}^{\zeta} ) \subset {\mathcal  V}_{\iota} \times {\mathcal V}_{\iota}  $, and such that
\begin{eqnarray}
 \chi_{\iota |_{\rm diagonal}}^{\zeta} = \chi_{\iota} . \label{additionparametrix}
\end{eqnarray}

\begin{defi} \label{pseudointro} 
For $ a^{\iota} \in S^m_{\iota} ({\mathcal V}_{\iota} \times \Ra^n) $, the pseudo-differential operator 
$$ \emph{op}_{w,h}^{\iota}(a^{\iota}) :
C_0^{\infty}(\Ra^n) \rightarrow
C_0^{\infty}({\mathcal V}_{\iota}) $$  is the operator  with
kernel
\begin{eqnarray}
 (2 \pi)^{-n} \int \! \! \int e^{i(r-r^{\prime})\rho + i
(\theta-\theta^{\prime})\cdot \eta} a^{\iota}(r,\theta,h \rho , h
w(r)\eta) d \rho d \eta \times
\chi_{\iota}^{\zeta}(r,\theta,r^{\prime},\theta^{\prime})
\label{supportpropreinfini}, \qquad \mbox{if} \ \iota \in I_{\infty}, \\
 (2 \pi)^{-n}  \int e^{i (x - x^{\prime} ) \cdot \xi } a^{\iota}(x,h \xi)
d \xi \times \chi_{\iota}^{\zeta}(x,x^{\prime}) , \qquad \qquad \qquad \mbox{if} \ \iota \in I_{\rm comp}.
\label{supportproprecompact}
\end{eqnarray}
\end{defi}

In other words, $
\mbox{op}_{h,w}^{\iota}(a^{\iota}) $ is obtained by multiplying the kernel of $ a^{\iota} (r,\theta,h D_r, h w(r)
D_{\theta}) \chi_{\iota} $ (resp. of $ a^{\iota}
(x,hD_x)  \chi_{\iota} $) by $ \zeta
((r,\theta)-(r^{\prime},\theta^{\prime})) $ (resp. by $ \zeta (x-x^{\prime}) $).   

If $ m < - n $ the integrals in (\ref{supportpropreinfini}) and (\ref{supportproprecompact})
are absolutely convergent, otherwise they must be understood as 
oscillatory integrals in the usual way.
That $ \mbox{op}^{\iota}_{w,h} (a_{\iota})
$ maps $ C_0^{\infty} (\Ra^n) $ into $ C_0^{\infty} ({\mathcal
V}_{\iota}) $ follows from the construction of $
\chi^{\zeta}_{\iota} $. Note also that
\begin{eqnarray}
 \mbox{op}^{\iota}_{w,h} (1)  = \chi_{\iota} , \label{identite}
\end{eqnarray}
 since the oscillatory integral is the Dirac measure along
the diagonal and $
\chi^{\zeta}_{\iota}(r,\theta,r^{\prime},\theta^{\prime}) =
\chi_{\iota}(r^{\prime},\theta^{\prime}) $ for $ |r-r^{\prime}| +
| \theta - \theta^{\prime} | $ small enough.

\medskip

\noindent {\bf Remark.} Note the factor $ w(r) $ in front of $ \eta $ in the amplitude of (\ref{supportpropreinfini}). 
The choice of notation of Definition \ref{pseudointro} is thus consistent with the expressions of the principal symbol of $ - \Delta_g $ given by (\ref{symbolprincipalinfini}) and (\ref{symbolprincipalcompact}).

\bigskip

We are now ready to state our results. We consider
$$ \varphi \in S^{-\sigma}(\Ra) , \qquad \sigma > 0  , $$
that is $ |\varphi^{(k)}(\lambda) | \leq C_k \scal{\lambda}^{-\sigma-k} $. The best constants $C_k$ are semi-norms defining the topology of $ S^{- \sigma}(\Ra) $.
\begin{theo} \label{theocalcfon1} Let $ P   $ denote
either $ - \Delta_g  $ or $ - \widetilde{\Delta}_g  $.  For all $ N
\geq 0$, the following holds:
$$ \varphi (h^2 P) = \sum_{\iota \in
\iota} {\mathcal Q}_N^{\iota} (P,\varphi,h) + h^{N+1} {\mathcal
R}_N (P,\varphi,h), \qquad h \in (0,1],
 $$
where, for all $ \iota \in I  $,
$$ \Psi_{\iota *} {\mathcal
Q}_{N}^{\iota}(P,\varphi,h) \Psi_{\iota}^* = \sum_{j=0}^N h^j
\emph{op}_{w,h}^{\iota}(a_j^{\iota})
$$
 with symbols $ a_0^{\iota} , \ldots , a_N^{\iota} $ of the form 
\begin{eqnarray}
a_0^{\iota} = \varphi \circ p_2^{\iota} , \qquad a_j^{\iota}  = \sum_{ k
\leq k (j)  \ } d_{jk}^{\iota} \varphi^{(k)} \circ p_2^{\iota} ,
\ \ \ j \geq 1 , \label{formeexplicitesymbole}
\end{eqnarray}
using the functions $
p_2^{\iota} $  given by (\ref{symbolprincipalinfini}) for $
\iota \in I_{\infty} $ and (\ref{symbolprincipalcompact})
for $ \iota \in I_{\rm comp} $. Here $ k (j) < \infty $
and
$$ d_{jk}^{\iota} \in S^{2k-j}_ {\iota}( {\mathcal V}_{\iota} \times \Ra^n) \
$$ 
is polynomial in the momentum variable ($
d_{jk}^{\iota} \equiv 0 $ if $ 2 k - j < 0 $) and independent of $ \varphi $.

 In addition, for all $ m , m^{\prime} \in \Na $,  all $
A \in \emph{Diff}_w^m ({\mathcal M}) $, $ B \in
\emph{Diff}_w^{m^{\prime}}({\mathcal M}) $,  all $ p \in
[2,\infty] $  and all  $  N $ such that $ N > n - 2 \sigma + m +
m^{\prime}  $, we have
\begin{eqnarray}
\left| \left|  h^m A {\mathcal R}_{N} (-\Delta_g,\varphi,h)
h^{m^{\prime}} B \right| \right|_{L^2({\mathcal M},dg)\rightarrow
L^{p}({\mathcal M},dg)} & \lesssim & h^{-n \left( \frac{1}{2} -
\frac{1}{p} \right) },  \label{Sobolevnaturel}
\end{eqnarray}
and, for $ P = - \widetilde{\Delta}_g $,
\begin{eqnarray}
\left| \left| w(r)^{  \frac{n-1}{2} - \frac{n-1}{p} } h^m A {\mathcal
R}_{N} (-\widetilde{\Delta}_g,\varphi,h) h^{m^{\prime}} B \right|
\right|_{L^2({\mathcal M},\widetilde{dg})\rightarrow
L^{p}({\mathcal M},\widetilde{dg})} & \lesssim & h^{-n \left(
\frac{1}{2} - \frac{1}{p} \right) },  \ \ \
\label{Sobolevnatureltilde}
\end{eqnarray}
both for $  h \in (0,1]  $.
\end{theo}

This theorem roughly means that, near infinity, $ \varphi (h^2 P) $ is well approximated by
pseudo-differential operators with symbols of the form $ a (r,\theta,\rho,w(r)\eta) $. The principal symbol is for instance
$$ \varphi (p_2^{\iota}(r,\theta,\rho,w(r)\eta))  . $$
Note that, when $ \varphi \in C_0^{\infty}(\Ra)  $,  this symbol is compactly supported with respect to $
\rho $ but not uniformly with respect to $  \eta$: if $ w(r) \rightarrow 0 $ as $ r \rightarrow
\infty $, $ \eta $ is not confined in a fixed compact set, since
we only have $ |\eta| \lesssim w (r)^{-1} $.

The estimates  (\ref{Sobolevnaturel})  and  (\ref{Sobolevnatureltilde})
 follow from the Sobolev embedding  $ D ( (- \Delta_g
)^{k} ) \subset L^{\infty}({\mathcal M}) $ for $ k > n / 4 $ (see Proposition \ref{injectionlinfini}) and,
to that extent, Theorem \ref{theocalcfon1} is an $ L^2 $ theorem. 

\bigskip

We now consider the $ L^p \rightarrow L^p $ properties. Recall first a classical definition
\begin{defi} \label{definitionpoidstempere} A function $ W : \Ra \rightarrow (0,+\infty) $ is a temperate weight if, for some positive constants $ C$,$M $,
\begin{eqnarray}
W (r^{\prime}) \leq C W (r) (1+|r-r^{\prime}|)^M  , \qquad r,r^{\prime} \in \Ra . \label{poidspolynomial}
\end{eqnarray}
\end{defi}
  The meaning of this definition is that $ W $ can neither grow nor decay too fast. For instance if $ d^k
w^{-1} / dr^k $ is bounded on $ \Ra $,  $w$ is a temperate weight.  This is an elementary consequence of Taylor's formula to order $k$ 
and of the fact that $ | d^j w^{-1} / dr^j | \lesssim
w^{-1}  $, by  (\ref{cinfiniborne}).

The operators $ \mbox{op}_{w,h}^{\iota}(a_j^{\iota}) $ of Theorem \ref{theocalcfon1}
are bounded on $ L^p ({\mathcal M},dg) $, $ L^p ({\mathcal M},\widetilde{dg}) $, or more generally on   $ L^p ({\mathcal M},W(r) dg) $ and $ L^p ({\mathcal M}, W(r) \widetilde{dg}) $ for all temperate weight $W$ and all $ p \in [1,\infty] $ (see Proposition \ref{appSchur}).
We  therefore focus on the remainder terms $ {\mathcal R}_N (P,\varphi,h) $.

\begin{theo} \label{Lprestetilde}  For all $ N \geq 0 $, all temperate weight $W$ and all $ 1 < p < \infty $,
\begin{eqnarray}
\left| \left| W(r)^{-1} {\mathcal R}_{N} (-\widetilde{\Delta}_g,\varphi,h) ) W(r)
\right| \right|_{L^p({\mathcal M},\widetilde{dg})\rightarrow
L^{p}({\mathcal M},\widetilde{dg})} \leq C_{N,p,\varphi,W}, \qquad h \in
(0,1] . \label{Lppourlerestetilde}
\end{eqnarray}
The constant $ C_{N,p,\varphi,W} $ depends (linearly) on a finite number of semi-norms of $ \varphi \in S^{- \sigma}(\Ra) $.
\end{theo}

\begin{coro} \label{coro1}   For all $ 1 < p < \infty $ and all temperate
  weight $W$,
$$ \left| \left|W(r)^{-1}  \varphi (-h^2\widetilde{\Delta}_g) W(r)
\right| \right|_{L^p({\mathcal M},\widetilde{dg})\rightarrow
L^{p}({\mathcal M},\widetilde{dg})}  \lesssim  1, \qquad h \in
(0,1] . $$
Equivalently, we have
$$ \left| \left| W(r)^{-1}  w(r)^{\frac{n-1}{p} - \frac{n-1}{2}}  \varphi
    (-h^2\Delta_g) w(r)^{\frac{n-1}{2} - \frac{n-1}{p}} W(r)
\right| \right|_{L^p({\mathcal M},dg)\rightarrow
L^{p}({\mathcal M},dg)}  \lesssim  1, \qquad h \in
(0,1] . $$
\end{coro}

Observe that Theorem \ref{Lprestetilde}  and Corollary 
\ref{coro1} hold in particular if $ w (r) = e^{-r} $ in which case $ \varphi (-h^2 \Delta_g) $ is in general not bounded on $ L^p ({\mathcal M},dg) $. 
Theorem \ref{Lprestetilde} is a consequence of a stronger  result, namely  Proposition 
\ref{impliqueLprestetilde}, showing  that, in any
chart, the resolvent $ (z-\widetilde{\Delta}_g)^{-1}  $ is a
pseudo-differential operators whose full symbol belongs to a suitable class. 
Since this result is of more technical nature, we prefer not to state it in this part.

\smallskip

If the function $w$ itself is a temperate weight, for instance if $ w (r) = r^{-1} $ for $r$ large, Theorem \ref{Lprestetilde} also implies the following result.

\begin{coro} \label{Lpreste} If $ w $ is a temperate weight, then for all temperate weight $ W $,
  all  $ N \geq 0 $  and all $ 1 < p < \infty $,
\begin{eqnarray}
\left| \left| W^{-1}(r) {\mathcal R}_{N} (-\Delta_g,\varphi,h) W(r) \right|
\right|_{L^p({\mathcal M},dg)\rightarrow L^{p}({\mathcal M},dg)} \leq C_{N,p,\varphi,W} , \qquad h \in (0,1] . \label{Lppourlereste}
\end{eqnarray}
The constant $ C_{N,p,\varphi,W} $ depends (linearly) on a finite number of semi-norms of $ \varphi \in S^{- \sigma}(\Ra) $.
In particular,
\begin{eqnarray}
 \left| \left| W(r)^{-1}  \varphi (-h^2\Delta_g) W(r) \right|
\right|_{L^p({\mathcal M},dg)\rightarrow L^{p}({\mathcal M},dg)}
\lesssim  1, \qquad h \in (0,1] . \label{explicitation}
\end{eqnarray}
\end{coro}
Of course, (\ref{explicitation}) holds with $ W = 1 $. As explained in the
introduction, this last result can be considered as essentially well known (see for instance \cite{Tayl0} for $ h=1 $). We 
quote it to emphasize the difference with Corollary \ref{coro1} where $ w $ is not assumed to be a temperate weight. 
It follows directly from Theorem \ref{Lprestetilde}, using (\ref{weightlp}),
(\ref{laplacientilde}) and the fact that products or real powers of temperate
weights are temperate weights.

\section{Parametrix of the resolvent and applications} \label{pseudodiffw}
\setcounter{equation}{0}

In the main part of this section, namely until (\ref{aadditionner}), we
work in coordinate patches $ {\mathcal U}_{\iota} $ of the form
(\ref{ouverts}) (ie with $ \iota \in I_{\infty} $).

\subsection{Elementary pseudo-differential calculus}
In this part, we give elementary composition formulas and the related remainder estimates for pseudo-differential operators of the form $ \mbox{op}_{w,h}^{\iota}(a) $.
We will note develop a systematic study of  the symbolic calculus
but only record the basic results required for the calculation of parametrices of
 $ ( z - h^2\Delta_g )^{-1} $ and $
(z- h^2\widetilde{\Delta}_g)^{-1} $.

For $ \Omega \subset \Ra^D $, $ D \geq 1 $,  $ C_b^{\infty}(\Omega) $ will denote the space of smooth
functions bounded on $ \Omega $ as well as their
derivatives.

For  $ b \in S^m_{\iota} ( {\mathcal V}_{\iota}
\times \Ra^n) $ and $ h \in (0,1] $, we  set
\begin{eqnarray}
 \left[ b (r,\theta,h D_r, h w(r)D_{\theta}) v \right] (r,\theta)
= (2\pi)^{-n} \int \! \! \int e^{ i (r\rho + \theta.\eta)} b
(r,\theta,h \rho , h w (r) \eta) \widehat{v}(\rho,\eta) d\rho
d\eta \label{pseudo}
\end{eqnarray}
with $ \hat{v}(\rho,\theta) = \int \! \! \int e^{- ir\rho - i
\theta.\eta} v (r,\theta) d r d \theta $ the usual Fourier
transform. %Observe that this  operator is well defined from $
%C_0^{\infty}(\Ra^n) $ to $ C^{\infty} (\Omega) $, since $ w $ is
%defined  on $ \Omega $. 
In the special case of a polynomial symbol
in $ \rho $ and $ \eta $, $  a (r,y,\rho,\eta) = \sum
a_{j\alpha}(r,\theta) \rho^j \eta^{\alpha} $, we have
\begin{eqnarray}
  a (r,\theta,h D_r, h w(r)D_{\theta}) = \sum
a_{j\alpha}(r,\theta) (h w(r) D_{\theta})^{\alpha} (h D_r)^j ,
\label{opdif}
\end{eqnarray}
where one must notice that $ D_r $ and $ w(r)D_{\theta} $ don't
commute.

 We have the following elementary
result.
\begin{prop} \label{composition} % Let
%$ \Omega \subset (R,+\infty) \times \Ra^{n-1} $ be an open set. 
Let $ a \in S^{m_1}_{\iota}( {\mathcal V}_{\iota} \times \Ra^n) $ be  polynomial
 in $ (\rho , \eta) $ and let $ b \in S^{m_2} ( {\mathcal V}_{\iota} \times
\Ra^n) $ with $ m_2 \in \Ra $. We have
\begin{eqnarray}
a (r,\theta,h D_r, h w(r)D_{\theta})  b (r,\theta,h D_r, h
w(r)D_{\theta}) = \sum_{l = 0}^{m_1} h^l (a \# b)_l (r,\theta,h D_r,
h w(r)D_{\theta}) \label{comp2}
\end{eqnarray}
where, if we set
$$ D_w = D_r + \frac{w^{\prime}(r)}{w(r)} \eta \cdot D_{\eta} , $$
the symbol $ (a \# b)_k = (a \# b)_k (r,\theta,\rho,\eta) \in
S^{m_1 + m_2 - k} ( {\mathcal V}_{\iota} \times \Ra^n) $ is given by
$$ (a \# b)_k = \sum_{j + |\beta| = k} \frac{1}{j! \beta !}  w(r)^{|\beta|}
\left( \partial_{\rho}^j \partial_{\eta}^{\beta} a \right) \left(
D_{\theta}^{\beta} D_w^{j} b \right) .
$$
\end{prop}
When $ w  \equiv 1 $, this proposition is of course the usual
composition formula for pseudo-differential operators. Note 
that, since $ a $ is polynomial of degree $ \leq m_1 $, we have $
(a \# b)_l \equiv 0 $ for $ l > m_1 $ and the composition formula is exact (there is no remainder term).

\bigskip

\noindent {\bf Remark.} 
A simple induction shows that the
operator $ D_w^{j} $ is a linear combination of
\begin{eqnarray}
 \left( \frac{w^{\prime}(r)}{w(r)} \right)^{(j_1)} \cdots \left(
 \frac{w^{\prime}(r)}{w(r)} \right)^{(j_k)} D_r^{l} \eta^{\alpha}
 D_{\eta}^{\alpha} \label{deriveetordue}
\end{eqnarray}
with  $  j_1 + \cdots + j_k + k + l = j  $, $ |\alpha| \leq k $
and $ k \geq 0 $. If $ k = 0 $ then $ (w^{\prime}/w)^{(j_1)}
\cdots (w^{\prime}/w)^{(j_k)} = 1 $. The notation $
(w^{\prime}/w)^{(j_i)} $ stands for the $ j_i $-th derivative of $
w^{\prime}/w$ with $ j_i \geq 0 $.

\bigskip

\noindent {\it Proof of Proposition \ref{composition}.}  Applying the right and side of (\ref{opdif}) to $
\reff{pseudo} $, the result  follows from the Leibniz rule and the fact that
$$ D_r \left( b
(r,\theta,h \rho , h w (r) \eta) \right) = \left( D_w b \right)
(r,\theta,h \rho , h w (r) \eta) .
$$
We omit the standard details of the calculation. That $
(a\# b)_k $ belongs to $ S^{m_1 + m_2 - k}( {\mathcal V}_{\iota} \times \Ra^n) $
follows from (\ref{cinfiniborne}) using (\ref{deriveetordue}). \finpreuve

\bigskip

We next consider the  pseudo-differential
quantization $ \mbox{op}_{w,h}^{\iota}(\cdot) $ given by (\ref{supportpropreinfini}).

\begin{prop} \label{complocale}
Let $ a \in S^{m_1}_{\iota}( {\mathcal V}_{\iota} \times \Ra^n) $ be  polynomial
 in $ (\rho , \eta) $ and let $ b \in S^{m_2} ( {\mathcal V}_{\iota} \times
\Ra^n) $ with $ m_2 \in \Ra $. Let $ W $ be a positive function on $ \Ra $ such that
\begin{eqnarray}
W(r) \leq C W (r^{\prime}), \qquad |r-r^{\prime}| \leq 1 . \label{superpoids}
\end{eqnarray}
Then, for all $ N
> 0 $,
$$  a (r,\theta,h D_r, h w(r)D_{\theta}) \emph{op}_{w,h}^{\iota}(b) = \sum_{l=0}^{m_1} h^{l} \emph{op}_{w,h}^{\iota}\left( (a\#b)_l
\right) + h^{N+1} \overline{R}_N^{\iota} (h,a,b),
$$
where, for all  $
k_1,k_2 \in \Na $, all $ A_1 \in \emph{Diff}_w^{k_1} ({\mathcal M}) $,  $ A_2 \in \emph{Diff}_w^{k_2} ({\mathcal M}) $ and all $ p \in [1,\infty ] $,
\begin{eqnarray}
\left| \left| W(r) A_1 \Psi_{\iota}^* \overline{R}_N^{\iota} (h,a,b) \Psi_{\iota *} A_2 W(r)^{-1} \right|
\right|_{L^{p}({\mathcal M},\widetilde{dg}) \rightarrow L^{p}({\mathcal M},\widetilde{dg})} 
%\leq  C_{p k
 % k^{\prime} \beta \beta^{\prime}}  ,
 \lesssim 1 , \label{estimereste} \\
 \left| \left| w(r)^{\frac{n-1}{2}} W(r) A_1 \Psi_{\iota}^* \overline{R}_N^{\iota} (h,a,b) \Psi_{\iota *} A_2 W(r)^{-1} \right|
\right|_{L^{2}({\mathcal M},\widetilde{dg}) \rightarrow L^{\infty}({\mathcal M})}
% \leq C^{\prime}_{p  k^{\prime} \beta \beta^{\prime}}
\lesssim 1 , \label{estimereste2}  
\end{eqnarray}
for $ h \in (0,1 ] $.
More precisely the norms in (\ref{estimereste}) and (\ref{estimereste2})
 are controlled by a finite number of
semi-norms of $ a $ and $ b $ independent of $h$.
\end{prop}

Note that the condition (\ref{superpoids}) is satisfied if $ W $ is a temperate weight but also by any power of $w$. In particular, $ W (r) = e^{\gamma r} $ is a possible choice although it is not a temperate weight. In particular, (\ref{estimereste}) and (\ref{estimereste2}) are respectively equivalent to 
\begin{eqnarray}
 \left| \left| W(r)  A_1  \Psi_{\iota}^* \overline{R}_N^{\iota} (h,a,b) \Psi_{\iota *}
 A_2  W(r)^{-1} \right|
\right|_{L^{p}({\mathcal M},dg) \rightarrow L^{p}({\mathcal M},dg)} 
%\leq  C_{p k
 % k^{\prime} \beta \beta^{\prime}}  ,
 \lesssim 1 , \label{estimerestebis} \\
 \left| \left|  W(r) A_1 \Psi_{\iota}^* \overline{R}_N^{\iota} (h,a,b) \Psi_{\iota *} A_2 W(r)^{-1} \right|
\right|_{L^{2}({\mathcal M},dg) \rightarrow L^{\infty}({\mathcal M})}
% \leq C^{\prime}_{p  k^{\prime} \beta \beta^{\prime}}
\lesssim 1 , \label{estimereste2bis}
\end{eqnarray}
They are simply obtained by replacing $ W(r) $ respectively by $ W (r) w(r)^{\frac{1-n}{p}} $ and $ W (r) w(r)^{\frac{1-n}{2}} $ which both satisfy (\ref{superpoids}).

By opposition to Proposition \ref{composition}, we now have 	a remainder. It is  due to the derivatives of cutoff near the diagonal in the definition of $ \mbox{op}_{w,h}^{\iota}(\cdot) $  but not to the tail of the expansion $
\sum_l h^l (a \# b)_l $ for this sum is finite.

Before proving this proposition, we state two lemmas which  will
be useful further on and whose proofs are very close to the proofs
of the estimates  (\ref{estimereste}) and  (\ref{estimereste2}).

\begin{lemm} \label{appSchur} Let $ c \in S^{m}_{\iota} ( {\mathcal V}_{\iota} \times \Ra^n) $ with $ m < 0
$ and let $ W $ be a positive function satisfying (\ref{superpoids}). 
 Then, for all $ p \in [1,\infty] $, we have
$$  \left| \left| W(r) \emph{op}_{w,h}^{\iota}\left( c \right)W(r)^{-1} \right|
\right|_{L^{p}(\Ra^n) \rightarrow L^{p}(\Ra^n)} \lesssim  1  ,
\qquad h \in (0,1] . $$
\end{lemm}

\bigskip

\noindent {\it Proof.} Consider first the case where $ W \equiv 1 $. If $ \hat{c} $ is the Fourier
transform of $c$ with respect to $ \rho,\eta $, the kernel of $
\mbox{op}_{w,h}^{\iota}\left( c \right) $ reads
$$  C_{\iota}(r,\theta,r^{\prime},\theta^{\prime},h) =
  h^{-n} w(r)^{1-n} \hat{c} \left( r,\theta, \frac{r^{\prime}-r}{h} ,
\frac{\theta^{\prime}-\theta}{h w(r)} \right) \frac{W(r)}{W(r^{\prime})}
\chi^{\zeta}_{\iota} (r,\theta,r^{\prime},\theta^{\prime}) .
$$
For $ (r,\theta) \in {\mathcal V}_{\iota} $,
$  c (r,\theta, .,.) \in L^{\epsilon + n / |m|} (\Ra^n_{\rho,\eta})   $,
with norm uniformly bounded with respect to $ (r,\theta)  $, thus $ \hat{c}
(r,\theta,.,.) $ belongs to a bounded subset of  $
L^1_{loc}(\Ra^n_{\hat{\rho},\hat{\eta}})  $
by Young's theorem.  Therefore, for all $ N $ we can write
\begin{eqnarray}
 | \hat{c}(r,\theta,  \hat{\rho} ,  \hat{\eta})| \leq C_N (1+f_0(r,\theta,\hat{\rho},\hat{\eta})) (
|\hat{\rho}| + |\hat{\eta}| + 1)^{-N} , \qquad
 (r,\theta) \in {\mathcal V}_{\iota}, \ \ \hat{\rho} \in \Ra, \ \hat{\eta} \in \Ra^{n-1}  ,
 \label{decayFouriertransform}
\end{eqnarray}
with $ f_0(r,\theta,.,.) $ bounded in $ L^1_{\rm comp}(\Ra^n_{\hat{\rho},\hat{\eta}}) $
Thus, the family  $ \hat{c}(r,\theta,.,.)  $ is bounded in $ L^1
(\Ra^n_{\hat{\rho},\hat{\eta}}) $. Elementary changes of variables show that
$$ \sup_{(r,\theta) \in \Ra^n} \int_{\Ra} \int_{\Ra^{n-1}} |C_{\iota}(r,\theta,r^{\prime},\theta^{\prime},h)|
 d r^{\prime}d \theta^{\prime} \lesssim 1 , \qquad
 \sup_{(r^{\prime},\theta^{\prime}) \in \Ra^n} \int_{\Ra}
 \int_{\Ra^{n-1}}  | C_{\iota}(r,\theta,r^{\prime},\theta^{\prime},h) | dr d\theta \lesssim 1 , $$
for $ h \in (0,1 ] $. Recall that $ C_{\iota} $ is globally
defined on $ \Ra^{2n} $ so the above quantities makes sense. The
result is then a consequence of the standard Schur lemma. For a general $W$ the same proof applies since we only have to multiply the kernel $ C_{\iota} $
by the bounded function $ W(r) \chi^{\zeta}_{\iota} (r,\theta,r^{\prime},\theta^{\prime}) W(r^{\prime})^{-1} $ on the support of which $ r -r^{\prime} $ is bounded. 
\finpreuve

\begin{lemm} \label{appSchurL2} Let $ c \in S^{m}_{\iota} ( {\mathcal V}_{\iota} \times \Ra^n) $ with $ m < - n/2
$ and let $ W $ be a positive function satisfying (\ref{superpoids}). Then
$$  \left| \left| w(r)^{\frac{n-1}{2}} W(r)  \emph{op}_{w,h}^{\iota}\left( c \right)W(r)^{-1} \right|
\right|_{L^{2}(\Ra^n) \rightarrow L^{\infty}(\Ra^n)} \lesssim  h^{-n/2}  ,
\qquad h \in (0,1] . $$
\end{lemm}

\noindent {\it Proof.} With the notation of the proof of Lemma \ref{appSchur},
 the result is a direct consequence of the estimate
\begin{eqnarray}
 \sup_{(r,\theta) \in \Ra^n} \int_{\Ra} \int_{\Ra^{n-1}} |w (r)^{\frac{n-1}{2}} W(r) C_{\iota}(r,\theta,r^{\prime},\theta^{\prime},h)
 W(r^{\prime})^{-1}|^2
 d r^{\prime}d \theta^{\prime} \lesssim h^{-n}, \qquad h \in (0,1] \nonumber 
\end{eqnarray}
which follows again from elementary changes of variables, using  that $ \hat{c}(r,\theta,.,.) $ belongs to a bounded subset of $ L^2 (\Ra^n) $
as $ (r,\theta) $ varies and that $ W (r)/W(r^{\prime}) $ is bounded on the support of $ C_{\iota} $. \finpreuve

\bigskip

\noindent {\bf Remark.} The proofs of both lemmas  still
hold if the kernel of $ \mbox{op}_{w,h}^{\iota}\left( c \right) $ is
multiplied by a bounded function. We shall use it in the
following proof.

\bigskip

\noindent {\it Proof of Proposition \ref{complocale}.}
We may clearly assume that (\ref{opdif}) is reduced to one
term.  Applying this operator to  (\ref{supportpropreinfini}) (with $a=b $)
 on the $r,\theta$
variables, we get  the kernel of $  \sum_k h^{k}
\mbox{op}_{w,h}^{\iota}\left( (a \# b)_k \right) $ (using Proposition
$\refe{composition}$) plus a linear combination of integrals of
the form
$$
a_{j\alpha} (r,\theta)  \int \! \! \! \int e^{i (r-r^{\prime})\rho
+ i(\theta - \theta^{\prime}).\eta} (h \rho)^{j_1} (h
\eta)^{\alpha_1}
 ( \partial^{\alpha_2}_{\theta} D_w^{j_2 } b)
(r,\theta,h \rho , h w(r)\eta) \ d\rho d \eta
 \partial_r^{j_3} \partial_{\theta}^{\alpha_3}
\chi^{\zeta}_{\iota}(r,\theta,r^{\prime},\theta^{\prime})  $$
where $ j_1 + j_2 +j_3 = j $, $ \alpha_1 + \alpha_2 + \alpha_3 =
\alpha  $ and $ j_3 + |\alpha_3| \geq 1 $. The latter implies that
$ \partial_r^{j_3} \partial_{\theta}^{\alpha_3}
\chi^{\zeta}_{\iota} $ is supported in  $
|(r,\theta)-(r^{\prime},\theta^{\prime})| \geq \varepsilon $ which
allows to integrate by parts using $
|(r,\theta)-(r^{\prime},\theta^{\prime})|^{-2} \Delta_{\rho,\eta}
  $. We thus obtain integrals of the form
\begin{eqnarray}
  h^{2N }\int \! \! \! \int e^{i (r-r^{\prime})\rho + i(\theta -
\theta^{\prime}).\eta}  c_N (r,\theta,h \rho , h w(r)\eta) \ d\rho
d \eta \frac{B_N
  (r,\theta,r^{\prime},\theta^{\prime})}{ |(r,\theta)-(r^{\prime},\theta^{\prime})|^{2N}
} \label{pseudohorsdiagonale}
\end{eqnarray}
 with $ N $ as large as we want, $ c_N \in S^{m+|\alpha| +j -
2 N} ( {\mathcal V}_{\iota} \times \Ra^n) $ and $ B_N \in
C_b^{\infty}(\Ra^{2n}) $  with support in   $ \{ \varepsilon \leq
|(r,\theta)-(r^{\prime},\theta^{\prime})| \leq 2  \varepsilon \}
$. 
With no loss of generality, we may assume that
$$  \Psi_{\iota *} A_1 \Psi_{\iota}^* = (w(r)D_{\theta})^{\beta} D_r^k , \qquad    \Psi_{\iota *} A_2 \Psi_{\iota}^* = 
(w(r)D_{\theta})^{\beta^{\prime}} D_r^{k^{\prime}} . $$
Applying $ (w(r)D_{\theta})^{\beta} D_r^k $ to 
(\ref{pseudohorsdiagonale}) yields an integral of the same
form, using the boundedness of $ w $ and its derivatives.
 To apply (the transpose of)  $
(w(r^{\prime})D_{\theta^{\prime}})^{\beta^{\prime}}
D_{r^{\prime}}^{k^{\prime}} $ to the kernel of $ \overline{R}_N^{\iota} (a,b,h) $, we rewrite this operator as $ (w (r^{\prime})
/ w(r) )^{|\beta^{\prime}|}
(w(r)D_{\theta^{\prime}})^{\beta^{\prime}}
D_{r^{\prime}}^{k^{\prime}}  $. We still obtain integrals
of the same form as (\ref{pseudohorsdiagonale}) multiplied by derivatives of  $  (  w (r^{\prime})/w (r)
)^{|\beta^{\prime}|} $. By  (\ref{diag}), these derivatives   are
bounded since $ | r-r^{\prime} | \leq 2 \varepsilon $ on the
support of $ B_N $. Then  (\ref{estimereste})  and 
(\ref{estimereste2})  follow respectively from the proofs of Lemma
 \ref{appSchur}  and \ref{appSchurL2}. \finpreuve

\bigskip

So far, we have considered composition with differential operators
to the left. Since our operators are properly supported, the
composition to the right can be also easily considered.

\begin{prop} \label{complocale2} Let $a$ and $b$ be as in Proposition  \ref{complocale} and let $W$ be a positive function satisfying (\ref{superpoids}).
Then, for all $ N > m_1 + m_2 + n  $, we have $$
\emph{op}_{w,h}^{\iota}(b) a(r,\theta,hDr,hw(r)D_{\theta}) = \sum_{l= 0 }^{N} h^{l} \emph{op}_{w,h}^{\iota}\left(
c_l \right) + h^{N+1} \underline{R}_N^{\iota} (h,a,b) $$ with $
c_l \in S^{m_1 + m_2 - l}_{\iota}({\mathcal V}_{\iota} \times \Ra^n) $ depending continuously on $a$ and $b$,
and $ \underline{R}_N^{\iota} (h,a,b) $ an operator with
continuous kernel supported in $ {\mathcal V}_{\iota} \times
{\mathcal V}_{\iota} $. Moreover, for all $N$,  all  $ k_1,k_2 \in \Na $ such
that
$$ N > m_1 + m_2 + n + k_1 + k_2  ,   $$
all $ A_1 \in \emph{Diff}_w^{k_1} ({\mathcal M}) $, $ A_2 \in \emph{Diff}_w^{k_2} ({\mathcal M}) $ and for all $ p \in [1,\infty ] $, we have
\begin{eqnarray}
 \left| \left|W(r) A_1 \Psi_{\iota}^* \underline{R}_N^{\iota} (h,a,b) \Psi_{\iota *} A_2 W(r)^{-1} \right|
\right|_{L^{p}({\mathcal M},\widetilde{dg}) \rightarrow L^{p}({\mathcal M},\widetilde{dg})} \lesssim 1
%\leq  C_{p k
%  k^{\prime} \beta \beta^{\prime}}  ,
%\qquad h \in (0,1] 
, \nonumber \\
 \left| \left| W(r) w(r)^{\frac{n-1}{2}} A_1 \Psi_{\iota}^* \underline{R}_N^{\iota} (h,a,b) \Psi_{\iota *}
A_2 W(r) \right|
\right|_{L^{2}({\mathcal M},\widetilde{dg}) \rightarrow L^{\infty}({\mathcal M})} \lesssim 1
% C^{\prime}_{p k
%  k^{\prime} \beta \beta^{\prime}}  h^{-n/2} ,
%\qquad h \in (0,1] 
, \nonumber
\end{eqnarray}
for $ h \in (0,1 ]$. More precisely, these norms are controlled by a finite number of semi-norms of $a$ and $b$ independent of $h$.
\end{prop}

 We will not need the explicit forms of the
symbols $c_l$ since we will only use this proposition for the analysis of some remainder terms.

Note also that the estimates on $ \underline{R}_N^{\iota} (h,a,b) $ have analogues with respect to the measure $dg$, similar to (\ref{estimerestebis}) and (\ref{estimereste2bis}),

\bigskip

\noindent {\it Proof.} We have to apply the transpose of $ a
(r^{\prime},\theta^{\prime},h D_{r^{\prime}},h
w(r^{\prime})D_{\theta^{\prime}}) $ to the Schwartz kernel of $
\mbox{op}_{w,h}^{\iota}(b) $. For simplicity we assume first that $ a
(r^{\prime},\theta^{\prime},\rho,\eta) = w(r^{\prime}) \eta_1 $.
 By Taylor's formula, we
have
$$ w (r^{\prime}) = w (r) \left( 1 + \sum_{j = 1}^{N} \frac{1}{j !}
\frac{w^{(j)}(r)}{w(r)} (r^{\prime}-r)^{j} +  \frac{(r^{\prime} -
r)^{N+1}}{N!} \int_0^1 (1-t)^{N} \frac{w^{(N+1)}(r+ t
(r^{\prime}-r))}{w(r)} dt \right) .
$$
Integrating by parts with respect to $ \rho $ in the kernel of $
\mbox{op}_{w,h}^{\iota}(b) $, the principal part of the Taylor expansion
yields the expected expansion with
$$ c_l (r,\theta,\rho,\eta) = \frac{1}{j!} D_{\rho}^{j} b (r,\theta,\rho,\eta)
\frac{w^{(j)}(r)}{w(r)} \eta_1 . $$ The remainder is given by two
types of terms:  first by the derivatives $ D_{\theta_1^{\prime}}
$ falling on $ \chi^{\zeta}(r,\theta,r^{\prime},\theta^{\prime})
$, which yields kernels of the form  (\ref{pseudohorsdiagonale}), 
and second by the remainder in the Taylor formula thanks to
which we can integrate by parts $ N $ times with respect to $ \rho
$. In this case, we get a kernel of the form (\ref{pseudohorsdiagonale}), with $ N $ instead of $ 2N $ and a
symbol $c_N \in S^{m_1 + m_2 -N}_{\iota} ({\mathcal V}_{\iota} \times
\Ra^n) $.
 Since $ r-r^{\prime} $ is bounded on the support of $
\chi^{\zeta} $, $ w^{(N)}(r+ t (r^{\prime}-r)) /w(r) $ is bounded
too, uniformly with respect to $t \in [0,1] $, and the study of
the remainder is similar to the one of Proposition \ref{complocale}. By induction, we obtain the result if $ a = (
w(r^{\prime}) \eta )^{\alpha} $. Derivatives with respect to $r$
or multiplication operators are more standard and  studied similarly. \finpreuve

\subsection{Parametrix of the resolvent}

In this subsection, we construct a parametrix of the semi-classical resolvent of an operator $ P \in \mbox{Diff}^2_w ({\mathcal M}) $. Recall that this means that $ P $ is a differential operator of order 2 such that, in any
chart at infinity,
\begin{eqnarray}
  \Psi_{\iota *} P \Psi_{\iota}^*  = \sum_{k =
0}^{2}  p_{2-k}^{\iota} (r,\theta, D_r,  w(r)D_{\theta}) \label{notationlocale}
\end{eqnarray}
with $ p_{2 - k}^{\iota} \in S^{2-k}_{\iota}({\mathcal V}_{\iota} \times \Ra^n)  $.

We assume that
\begin{eqnarray}
P \ \mbox{is locally elliptic} , \label{cond0}
\end{eqnarray}
ie, in any chart, its principal symbol $ p^{\iota}_{\rm pr}
(x,\xi) $ satisfies $ | p_{\rm pr}^{\iota} (x,\xi) | \gtrsim
|\xi|^2 $ for $ \xi \in \Ra^n $, locally uniformly with respect to
$ x $. If $ \iota \in I_{\infty} $ , using the notation (\ref{notationlocale}), we furthermore assume that, for all $ K_{\iota} \Subset V_{\iota} $ (see (\ref{ouverts})),
\begin{eqnarray}
 | p_2^{\iota} (r,\theta,\rho, \eta) | \gtrsim \rho^2 + |\eta|^2,
 \qquad r > R, \ \theta \in K_{\iota}, \ \rho \in \Ra, \
 \eta \in \Ra^{n-1} . \label{elliptictordu}
\end{eqnarray}
Note that this is not a  lower bound  for
the principal symbol of $ \Psi_{\iota *} P \Psi_{\iota}^* $,
namely $ p^{\iota}_{2} (r,\theta, \rho,w(r)\eta) $, whose modulus
is only bounded from below by $ \rho^2 + w(r)^2 |\eta|^2 $. This is nevertheless the natural (degenerate) global ellipticity condition in this context.
We next
define $ {\mathcal C} \subset \Ca $ as
\begin{eqnarray}
{\mathcal C} = \mbox{closure of the range of the principal symbol
of } P \label{clos}
\end{eqnarray}
which is invariantly defined for the principal symbol is a
function on $ T^* {\mathcal M} $. We assume that $ {\mathcal C}
\ne \Ca $. In the final applications, with $ P = - \Delta $ or $ - \widetilde{\Delta}_g $, we will of course have $ {\mathcal C} = [ 0 , + \infty ) $.

We now seek an approximate inverse of $ h^2 P - z $, for $
h \in (0,1 ] $   and  $ z \in \Ca \setminus {\mathcal C} $.

\medskip

We work first in a patch at
infinity. Using the notation of (\ref{notationlocale}), we set
for simplicity
$$ p_{2} = p_2^{\iota} - z , \qquad p_{1} = p_1^{\iota} , \qquad p_{0} = p_0^{\iota} . $$
Observe that $ p_0,p_1 $ don't depend on $z$ but that $p_2 $ does.
We then have
$$ h^2  \Psi_{\iota *} P \Psi_{\iota}^* - z = \sum_{k =
0}^{2} h^k p_{2-k} (r,\theta,h D_r, h w(r)D_{\theta}) . $$
 For a given $ N \geq 0 $, we look for symbols $ q_{-2}, q_{-3} ,
\ldots , q_{-2-N} $ satisfying
\begin{eqnarray}
 \left( \sum_{k = 0}^{2} h^k p_{2-k} (r,\theta,h D_r, h
w(r)D_{\theta})  \right) \left( \sum_{j = 0}^{N} h^j
\mbox{op}_{w,h}^{\iota}( q_{-2-j} ) \right) = \chi_{\iota} + {\mathcal
O}(h^{N+1}) , \label{parametrixlocale}
\end{eqnarray}
where $ \chi_{\iota} $ is defined by (\ref{notationcutoff2}) and where  $ {\mathcal O} (h^{N+1}) $ will be given a precise
meaning below. Of course, we need to find such a family of
symbols for each patch, ie $ q_{-2-j} $ depends on $ \iota $, but
we omit this dependence for notational simplicity. By
Proposition \ref{complocale}, the left hand side of
(\ref{parametrixlocale}) reads
$$ \sum_{k+j + l \leq N} h^{k+j+l} \mbox{op}_{w,h}^{\iota} \left( ( p_{2 - k} \# q_{-2-j} )_l \right) +
h^{N+1}  R^{\iota}_N (h,z)  $$ 
where 
\begin{eqnarray}
  R^{\iota}_N (h,z) = \sum_{  k  + j + l \geq N+1}
h^{k+ j+l - N -1 } \mbox{op}_{w,h}^{\iota} \left( ( p_{2 - k} \# q_{-2-j} )_l
\right) + \sum_{k,j} \overline{R}_N^{\iota} (h, h^{k} p_{2-k},
h^j q_{-2-j}) , \label{restelocal}
\end{eqnarray}
with $ \overline{R}_N^{\iota} $ defined in Proposition 
\ref{complocale}. In the above sums, we have $ 0 \leq k \leq 2
$, $ 0 \leq j \leq N $ and $ 0 \leq l \leq 2 $. Thus, by
(\ref{identite}), requiring (\ref{parametrixlocale}) 
leads to the following equations for $ q_{-2} , \ldots , q_{-2-N}
$
$$ \sum_{k+l+j = \nu } ( p_{2 - k} \# q_{-2-j} )_l =
  \begin{cases}
    1 & \text{if} \ \nu  = 0, \\
    0 & \text{if} \ \nu \geq 1,
  \end{cases} \qquad 0 \leq \nu \leq N .
 $$
This system is triangular and, since $ (a \# b)_0 =
ab $, its unique solution is given recursively by
$$ q_{-2} = \frac{1}{p_2}, \qquad \qquad
q_{-2-j} = - \frac{1}{p_2} \sum_{k+j_1+l=j \atop j_1 < j} (p_{2-k}
\# q_{-2-j_1})_l \qquad \mbox{for} \ j \geq 1 . $$
\begin{prop} \label{formpara} For all $ j \geq 1 $, $ q_{-2-j} $ is  a
finite sum (with a number of terms $k(j)$ depending on $j$ but not
on $z$) of the form
$$ q_{-2-j} = \sum_{k = 1}^{k(j)}  \frac{d_{jk}}{p_2^{1+k}} $$
where, for each $k$, $ d_{jk} \in S^{2 k - j}_{\iota} ({\mathcal
V}_{\iota} \times \Ra^n) $ is a polynomial in $ \rho,\eta $ which
is independent of $ z $ (in particular $ d_{jk} \equiv 0 $ when $
2 k - j < 0 $). More precisely, the coefficients of these
polynomials are linear combinations of products of derivatives of
$ w $, $ w^{\prime}/w $ and of the coefficients of $ p_0 , p_1 $
and $
\partial^{\alpha} p_2 $ with $ \alpha \ne 0 $.
\end{prop}

\medskip

\noindent {\it Proof.} This follows from an induction using 
(\ref{deriveetordue}) and the fact that, for any multi-index $
\alpha \ne 0 $,  $ \partial^{\alpha} ( 1/p_2^{1+k} ) $ is a linear
combination of
$$ \frac{\partial^{\alpha_1} p_2 \cdots \partial^{\alpha_{k^{\prime}} } p_2 }{p_2^{1+k+k^{\prime}}},
$$
with $ \alpha_1 + \cdots + \alpha_{k^{\prime}} = \alpha $, $ 1
\leq k^{\prime} \leq |\alpha| $ and  $ \alpha_i \ne 0 $ for all $
i \in \{ 1 , \ldots , k
^{\prime} \} $. \finpreuve

\bigskip

 With the
notation  (\ref{restelocal}), we set
$$  {\mathcal R}^{\iota}_N (h,z)  =  \Psi_{\iota}^*
R_N^{\iota} (h,z) \Psi_{\iota *} . $$
\begin{lemm} \label{lemmetheoreme} Let $ d \mu $ denote either $ d g $ or $ \widetilde{d g}
$. Then, for all positive function $ W $ satisfying (\ref{superpoids}), all $ p \in [1,\infty] $ and all  $ N \geq 0
$, there exists $ \nu > 0 $ such that, for all $ A \in
\emph{Diff}_w^{m} ({\mathcal M}) $ and $ B \in
\emph{Diff}_w^{m^{\prime}} ({\mathcal M}) $ with  $ m + m^{\prime} - N < 0  $, we have
$$  \left| \left| W(r) h^m A  {\mathcal R}^{\iota}_N (h,z)
h^{m^{\prime}}B W(r)^{-1} \right| \right|_{L^p ({\mathcal M},d\mu)
\rightarrow L^p({\mathcal M },d\mu)} \lesssim \left(
\frac{1+|z|}{{\rm dist}(z, {\mathcal C} )} \right)^{\nu} , $$ for all $ h
\in (0,1] $ and all $ z \notin {\mathcal C} $.
\end{lemm}

\noindent {\it Proof.} We first assume that  $ A = B = 1 $ (and
that $ m = m^{\prime} = 0 $). By (\ref{localisationnoyauoperateur}), the kernel of $ R_N^{\iota}(h,z) $ is supported in $ ( (R,+\infty) \times K_{\iota} )^2 $
for some $ K_{\iota} \Subset V_{\iota} $. Thus, using the equivalence of norms (\ref{equivalenceLpMRn}),  
the result, with $ d\mu = \widetilde{dg} $, is a direct
consequence of the bound
\begin{eqnarray}
  \left| \left|  W(r) R^{\iota}_N (h,z) W(r)^{-1} \right| \right|_{L^p
(\Ra^n) \rightarrow L^p(\Ra^n)}
\lesssim \left( \frac{1+|z|}{{\rm dist}(z, {\mathcal C} )}
\right)^{\nu} , \qquad h \in (0,1], \ \ z \notin {\mathcal C},
\label{bound0}
\end{eqnarray}
 which
follows from  Proposition \ref{complocale}  and Lemma 
\ref{appSchur}   once noticed  that each semi-norm of $ q_{-2-j}
$ in $ S^{-2-j}_{\iota} ({\mathcal V}_{\iota} \times \Ra^n) $ is bounded
by some power of $ (1+|z|)/ {\rm dist}(z,{\mathcal C}) $. The
latter is due to Proposition  \ref{formpara}  and
$$  \left| \frac{1 + \rho^2 + \eta^2}{p_{2}^{\iota} - z} \right| \thickapprox
\left| \frac{1 + p_2^{\iota}}{p_2^{\iota}-z} \right| \lesssim
\frac{1+|z|}{{\rm dist}(z, {\mathcal C} )} ,
$$
in which we used (\ref{elliptictordu}).
When $ d \mu = dg $, we use the equivalence (\ref{equivalenceLpapoids})
so that it is now sufficient to get the bound  (\ref{bound0}) 
with $ R^{\iota}_N (h,z) $ replaced by $ w(r)^{\frac{1-n}{p}}
R^{\iota}_N (h,z) w(r)^{\frac{n-1}{p}} $. The latter is clear for
this amounts to multiply the kernel of $ R^{\iota}_N (h,z) $ by $
(w(r^{\prime})/w(r))^{(n-1)/p} $ (which is bounded, using the boundedness of $ r-r^{\prime} $  on the
support of $ \chi^{\zeta}_{\iota} $ and  (\ref{diag})) so the (proofs of)
Proposition \ref{complocale} and Lemma \ref{appSchur}
still hold.

For general $ A $ and $ B$, we use Propositions 
\ref{complocale}  and  \ref{complocale2}  so that we are
reduced to the previous case  with an operator of the same form as
$ R^{\iota}_N (h,z) $ except that the symbols of the first sum in
 (\ref{restelocal}) now belong to $ S^{ -N + m+m^{\prime}}_{\iota}
({\mathcal V}_{\iota} \times \Ra^n) $. We can apply Lemma $
\refe{appSchur} $ to this term and the result follows. \finpreuve
\bigskip

 Let us now define
\begin{eqnarray}
 Q^{\iota}_N (h,z)  =   \sum_{j =
0}^{N} h^j \mbox{op}_{w,h}^{\iota}( q_{-2-j} ), \qquad  {\mathcal
Q}^{\iota}_N (h,z)  =  \Psi_{\iota}^* Q_N^{\iota} (h,z)
\Psi_{\iota *} .  \nonumber
\end{eqnarray}
Then, with $ f_{\iota} $ given by (\ref{partition}), we
obtain the relation
\begin{eqnarray}
(h^2 P - z) {\mathcal Q}^{\iota}_N (h,z) = f_{\iota} + h^{N+1}
{\mathcal R}^{\iota}_N (h,z) . \label{aadditionner}
\end{eqnarray}
So far, we have always assumed that $ \iota \in I_{\infty}
$, ie worked in patches at infinity, but the same analysis
still holds for relatively compact patches, ie for $ \iota \in
I_{\rm comp} $. We don't give the details of the construction
in the latter case for two reasons: the first is that this is
essentially well known for this is like working on a compact
manifold and the second is that the proofs are formally  the same
with the simpler assumptions that $ w \equiv 1 $ and that $
\chi_{\iota} $ is compactly supported.

Thus, by setting
$$ {\mathcal Q}_N (h,z) =  \sum_{\iota \in I} {\mathcal Q}^{\iota}_N (h,z) , \qquad
{\mathcal R}_{N}(h,z) =  \sum_{\iota \in I} {\mathcal
R}^{\iota}_N (h,z) , $$ then summing the equalities 
(\ref{aadditionner}) over $ I $ and using  (\ref{partition}), Lemma \ref{lemmetheoreme} gives the following
result where we recall that $ {\mathcal C} $ is defined by (\ref{clos}).
\begin{theo} \label{parametrixegenerale} Let $ P \in \emph{Diff}_w^2 ({\mathcal M}) $ be a second order differential operator
 satisfying (\ref{cond0})  and  (\ref{elliptictordu}). Then, for all $ N
\geq 0 $,  we have
\begin{eqnarray}
 (h^2 P - z) {\mathcal Q}_N (h,z) = 1 + h^{N+1} {\mathcal R}_{N}(h,z)
, \qquad  h \in (0,1 ] , \ \ z \notin {\mathcal C} . \label{formaladjoint} 
\end{eqnarray}
  If $
d\mu $ denotes either $ dg $ or $ \widetilde{dg} $, and $ m ,
m^{\prime} \in \Na $ satisfy $ m + m^{\prime} < N  $, then
for all $ p \in [1,\infty] $ and for all positive function $ W $ satisfying (\ref{superpoids}), there exists $ \nu \geq 0 $ such that,
for  all $ A \in \emph{Diff}_w^m ({\mathcal M}) $ and $ B \in
\emph{Diff}_w^{m^{\prime}}({\mathcal M}) $, we have
\begin{eqnarray} \left| \left|W(r) h^m A  {\mathcal R}_N (h,z)
h^{m^{\prime}}B W(r)^{-1} \right| \right|_{L^p ({\mathcal M},d\mu)
\rightarrow L^p({\mathcal M },d\mu)} \lesssim \left(
\frac{1+|z|}{{\rm dist}(z, {\mathcal C} )} \right)^{\nu} , \label{estimeeinjective}
\end{eqnarray}
for all $ h\in (0,1] $ and all $ z \notin {\mathcal C} $.
\end{theo}

This theorem gives a parametrix of the  resolvent of $h^2 P$ under the natural  ellipticity conditions (\ref{cond0}) and (\ref{elliptictordu})
(recall that if $w$ is not bounded from below, this corresponds to a degenerate ellipticity).

From now on, we assume that

 $$ P \ \mbox{is self-adjoint with respect to} \ d \mu = d g \ \mbox{or} \ \widetilde{dg} . $$
This condition is actually equivalent to the symmetry of $ P $ on $ C_0^{\infty} ({\mathcal M}) $. Indeed, (\ref{formaladjoint}) and (\ref{estimeeinjective}) implies that $ h^2 P \pm i $ is injective  for $h$ small enough, which shows that $P$ is essentially self-adjoint.

 The resolvent $ (h^2 P -z )^{-1} $ is then well defined for
all $ z \notin \Ra $ and
\begin{eqnarray}
 (h^2 P - z)^{-1} = {\mathcal Q}_N (z,h) - h^{N+1} (h^2 P -
z)^{-1} {\mathcal R}_{N} (h,z), \qquad z \notin \Ra, \ \ h \in
(0,1] . \label{derivableenz}
\end{eqnarray}
 Theorem  \ref{parametrixegenerale}  implies, for $ z =
i $, that in the operator norm on $ L^2 ({\mathcal M},d \mu) $, we
have
$$ (h^2 P - i)^{-1} \left( 1 + {\mathcal O}(h^{N+1}) \right) = {\mathcal Q}_N (h,i) $$
and thus, for some $ h_0 > 0 $ small enough and some bounded
operator $ B_1 $ on $ L^2 ({\mathcal M},d\mu) $, we get
$$ (h_0^2 P - i)^{-1} = {\mathcal Q}_N (i,h_0) B_1 . $$
More generally, for $ k \geq 1 $, we can write
$$ (h^2 P -z)^{-k} = \frac{1}{(k-1)!}
\partial_z^{k-1} (h^2 P -z)^{-1} $$
so applying $ (k-1)!^{-1}
\partial_z^{k-1} $ to  (\ref{derivableenz})  shows that $ (h^2 P
- z)^{-k} $ reads
\begin{eqnarray}
 (k-1)!^{-1} \partial_z^{k-1} {\mathcal Q}_N (z,h)
+ h^{N+1} (h^2 P - z)^{-k} \sum_{j=0}^{k-1} \frac{1}{j!} (h^2 P -
z)^{j} \partial_z^j {\mathcal R}_N (z,h) , \label{developpementquisameliore}
\end{eqnarray}
 using the holomorphy
of $ {\mathcal Q}_N (z,h) $ and $ {\mathcal R}_N (z,h) $ with
respect to $ z \in \Ca \setminus \Ra $ which standardly follows
from Proposition $  \refe{formpara} $. Therefore, by choosing $N$
large enough so that the sum above is bounded on $ L^2  $ (uniformly in $h$) and choosing then $
h = h_0 $ small enough, we obtain
\begin{eqnarray}
 (h^2_0 P - i)^{-k} = (k-1)!^{-1} \partial_z^{k-1} {\mathcal Q}_N
(z,h_0)_{|z=i} B_k , \label{astucepuissance}
\end{eqnarray}
 for some operator $ B_k $ bounded on $ L^2
({\mathcal M},d \mu) $.

\begin{lemm} \label{doubleborneL2} For all $ A \in \emph{Diff}_w^{2k} ({\mathcal M}) $, $ A
\partial_z^{k-1}{\mathcal Q}_N (z,h_0)
$ is bounded on $ L^2 ({\mathcal M},dg) $ and $ L^2 ({\mathcal M},
\widetilde{dg} ) $.
\end{lemm}

\noindent {\it Proof.} Consider first the case of $ \widetilde{dg}
$. By Proposition \ref{formpara}, for all $ \iota \in I
$, $ Q_{\iota} :=
\partial_z^{k-1} {\mathcal Q}_N^{\iota}(z,h_0)   $ is of the form $ \Psi_{\iota}^*
\mbox{op}_{w,h}^{\iota}(q_{\iota})\Psi_{\iota *}
 $ for some symbol $ q_{\iota} \in S^{-2 k}_{\iota} ({\mathcal V}_{\iota} \times \Ra^n)
$. A direct calculation shows that $ \Psi_{\iota *} A
\Psi_{\iota }^* \mbox{op}_{w,h}^{\iota}(q_{\iota}) $ has a kernel of the
form
$$ (2 \pi)^{-n} \int e^{i (x-y)\cdot \xi} a^{\iota} (x,y,\xi) d \xi , $$
with $ a^{\iota} \in C^{\infty}_b (\Ra^{3n}) $. Hence, the
corresponding operator is
 bounded on $ L^2 (\Ra^n) $ by the Calder\`on-Vaillancourt
theorem and thus its pullback $ A Q_{\iota} $ is bounded on $ L^2
({\mathcal M},\widetilde{dg}) $. The boundedness of $ A Q_{\iota}
$ on $ L^2 ({\mathcal M},dg) $ is equivalent to the one of $
w(r)^{(1-n)/2} A Q_{\iota} w^(r){(n-1)/2} $ on $ L^2({\mathcal
M},\widetilde{dg}) $. The latter follows from the same reasoning
using since $ w(r)^{(1-n)/2} A w(r)^{(n-1)/2} \in
\mbox{Diff}^{2k}_w ({\mathcal M}) $ and $ Q_{\iota}
$ is properly supported. \finpreuve

\medskip
Using  (\ref{astucepuissance}), Lemma \ref{doubleborneL2}
and the Spectral Theorem, we obtain the following
result.
\begin{prop} \label{regularite2} Let $ P \in \emph{Diff}_w^2({\mathcal M}) $, satisfying $ \reff{cond0} $ and $
\reff{elliptictordu} $, be self-adjoint with respect to $ d \mu =
dg $ or $ \widetilde{dg} $ . Then, for all $ k \geq 1 $ and all $
A \in \emph{Diff}_w^{2k} ({\mathcal M}) $, we have
$$ \left| \left| h^{2k} A (h^2 P - z)^{-k} \right| \right|_{L^2 ({\mathcal M},d\mu) \rightarrow
L^2 ({\mathcal M},d\mu)} \lesssim \frac{\scal{z}^k}{| \emph{Im} \
z|^k}, \qquad z \notin \Ra, \ \ h \in (0,1 ] . $$
\end{prop}

\bigskip

In the same spirit, we will prove the following  Sobolev injections.

\begin{prop} \label{injectionlinfini} Let $ P $ be as in Proposition  \ref{regularite2}
 and let $ k > n/4 $ be an integer. Then, if $ P $ is
self-adjoint with respect to $ d \mu = d g $, we have
$$  || (h^2 P - z)^{- k} ||_{L^2 ({\mathcal M},dg)\rightarrow L^{\infty}({\mathcal M})}
\lesssim h^{-\frac{n}{2}} \frac{\scal{z}^k}{| \emph{Im} \ z|^k}, \qquad z
\notin \Ra, \ \ h \in (0,1 ] .
 $$
If it is self-adjoint with respect to $ d \mu = \widetilde{dg} $,
we have
$$ || w(r)^{\frac{n-1}{2}} (h^2 P - z)^{-k} ||_{L^2 ({\mathcal M},\widetilde{dg})\rightarrow
L^{\infty}({\mathcal M})} \lesssim h^{-\frac{n}{2}}\frac{\scal{z}^k}{|
\emph{Im} \ z|^k}, \qquad z \notin \Ra, \ \ h \in (0,1 ] . $$
\end{prop}

Of course, by taking the adjoints, we have the corresponding $
L^1 \rightarrow L^2 $ inequalities.

\bigskip

\noindent {\it Proof.} We assume that $ d \mu = d g $. By (\ref{astucepuissance}), Lemma \ref{appSchurL2} with $ W (r) = w(r)^{\frac{1-n}{2}} $ and the equivalence of norms (\ref{equivalenceLpapoids}),  $ (h^2_0 P -i)^{-k} $
is bounded from $ L^2 ({\mathcal M},dg) $ to $ L^{\infty}({\mathcal M}) $. Therefore, using the spectral theorem,
$$ || (h^2 P - i)^{- k} ||_{L^2 ({\mathcal M},dg)\rightarrow L^{\infty}({\mathcal M})} \leq C \left|\left| \frac{(h_0^2 P - i)^k}{(h^2 P - i)^k} \right| \right|
_{L^2 ({\mathcal M},dg)\rightarrow L^{\infty}({\mathcal M})}
\lesssim h^{-k}  . $$
 Using \ref{developpementquisameliore} with $ z = i $ and $ N $ large enough, the estimate above improves to
 $$ || (h^2 P - i)^{- k} ||_{L^2 ({\mathcal M},dg)\rightarrow L^{\infty}({\mathcal M})}
\lesssim h^{- n / 2}  , $$
using again Lemma \ref{appSchurL2} for the principal part of the expansion. We then obtain the result from the estimate
$$ \left|\left| \frac{(h_0^2 P - i)^k}{(h^2 P - z)^k} \right| \right|
_{L^2 ({\mathcal M},dg)\rightarrow L^{2}({\mathcal M},dg)}
\lesssim \frac{\scal{z}^k}{|\mbox{Im}(z)|^k}  . $$
The case of $ d \mu = \widetilde{dg} $ is similar. \finpreuve

\subsection{Proof of Theorem \ref{theocalcfon1}}
  We shall use the classical
Helffer-Sj\"ostrand formula
\begin{eqnarray}
  \varphi (H) = \frac{1}{ \pi} \int \! \! \! \int_{\Ra^2}
\bar{\partial} \widetilde{\varphi} (x+iy)
 (H- x - i y)^{-1}  d x d y \label{HelfferSjostrand}
\end{eqnarray}
with $ \bar{\partial} = (\partial_x + i \partial_y)/2$, valid for any self-adjoint operator $H$. Here $
\widetilde{\varphi} \in C^{\infty}(\Ca) $ is an almost analytic
extension of $ \varphi $, ie such that $
\widetilde{\varphi}|_{\Ra} = \varphi $ and $
 \bar{\partial} \widetilde{\varphi}(z) $ vanishes to sufficiently high order
 on the real axis.

A justification of this formula for $ \varphi \in
C_0^{\infty}(\Ra)  $ can be found in \cite{DiSj1}.  It is shown in
\cite{Davi1} that, if $ \varphi \in S^{-\sigma}(\Ra) $ with $
\sigma > 0 $,  (\ref{HelfferSjostrand}) holds with $
\widetilde{\varphi}_M $ defined by
\begin{eqnarray}
\widetilde{\varphi}_M (x+iy) & = &  \chi_0 (y/\scal{x})
\sum_{k=0}^M f^{(k)} (x)  \frac{(iy)^k}{k!}, \label{holomorphe}
\end{eqnarray}
with $ M \geq 1 $ and $ \chi_0 \in C_0^{\infty}(\Ra)  $ such that
$ \chi_0 \equiv  1 $ near $ 0 $. With this choice, one has
\begin{eqnarray}
 |\bar{\partial} \widetilde{\varphi}_M (x+iy) | \lesssim |y|^M /
\scal{x}^{\sigma + 1 + M } ,  \qquad x,y \in \Ra . \label{dbarre}
\end{eqnarray}
This implies in particular that, for all integers $ \nu_1 \geq 1 $,
$ \nu_2 \geq 0 $ and $ M  \geq \nu_1 + \nu_2 $, we have
\begin{eqnarray}
\int \! \! \int_{\Ra^2} | \bar{\partial} \widetilde{\varphi}_M
(x+iy) | \times |y|^{- \nu_1} \left( \frac{ 1+|x|+ |y| }{|y|}
\right)^{\nu_2} dx dy < \infty, \label{justificationGreen}
\end{eqnarray}
which is easily seen by splitting the integral into two parts,
where $ |y| \leq 1 $ or $ |y| > 1 $, using the fact that $ |y| /
\scal{x} $ is bounded on the support of $ \widetilde{\varphi}_M $
in the latter case. If $ \sigma  > 1  $ and $ M \geq \nu  $, we also have
\begin{eqnarray}
\int \! \! \int_{\Ra^2} | \bar{\partial} \widetilde{\varphi}_M
(x+iy) |  \left( \frac{ 1+|x|+ |y| }{|y|}
\right)^{\nu} dx dy < \infty . \label{justificationGreenbis}
\end{eqnarray}

\bigskip

\noindent {\it Proof of Theorem \ref{theocalcfon1}.}
Let $ \iota \in I $. The form of $ \Psi_{\iota *} {\mathcal
Q}_N^{\iota} (P,\varphi,h) \Psi_{\iota}^* $, namely $
\reff{formeexplicitesymbole} $, simply follows by plugging the
expansion $ \reff{derivableenz} $ into $ \reff{HelfferSjostrand} $
and applying Green's formula. For the latter we use Proposition $
\refe{formpara} $ (recalling that $ p_2 \equiv p_2^{\iota} - z
$). All the integrals makes sense by $ \reff{justificationGreen} $
if we choose $ \widetilde{\varphi}_M $ with $ M \geq \max_{j \leq
N } (k(j)+1) $.

 Let us now prove $ \reff{Sobolevnaturel} $ and $
\reff{Sobolevnatureltilde} $. Since the proofs are very similar we
only show $ \reff{Sobolevnatureltilde} $ and thus consider $ P = -
\widetilde{\Delta}_g  $.  Fix $ N \geq 0 $. For $ N^{\prime}
> N $ and $ M  $ large enough, both to be chosen later, we set
\begin{eqnarray}
 {\mathcal R}_{N^{\prime}} (P,\varphi,h)  = \frac{1}{ \pi} \int \! \! \! \int_{\Ra^2}
\partial_{\bar{z}} \widetilde{\varphi}_M (x+iy)
 (h^2 P - x-iy)^{-1} {\mathcal R}_{N^{\prime}} (x+iy,h)  d x d y . \label{restedecale}
\end{eqnarray}
We next fix two integers $ k > n/4 $, $  \tilde{m} \geq m / 2 $,
and rewrite $ h^m A (h^2 P - z)^{-1} $ (with $z=x+iy$) as
\begin{eqnarray}
  (h^2 P - i)^{-k} \left\{
 (h^2 P - i)^k h^m A (h^2 P - i)^{-k-\tilde{m}} \right\}
 (h^2 P - z)^{-1}  (h^2 P - i)^{k+\tilde{m}}   . \label{integrabiliteverticale}
\end{eqnarray}
Using
 Proposition $ \refe{regularite2} $, the term $ \{ \cdots \} $  is bounded on
 $ L^2 ({\mathcal M},\widetilde{dg}) $ uniformly with respect to $h$. 
 On the other hand, by Theorem $ \refe{parametrixegenerale} $, there exists
 $ \nu_2 > 0 $ such that
$$  \left| \left| (h^2 P - i)^{k+\tilde{m}}{\mathcal R}_{N^{\prime}} (z,h) h^{m^{\prime}}
B  \right| \right|_{L^2 ({\mathcal M},\widetilde{dg})\rightarrow
L^2 ({\mathcal M},\widetilde{dg}) } \lesssim \scal{z}^{\nu_2}/
|\mbox{Im} \ z|^{\nu_2} ,
  $$
for $ z \notin \Ra $ and $ h \in (0,1] $,
  provided
\begin{eqnarray}
  N^{\prime}
> m^{\prime} + 2 (k + \tilde{m}) . \label{choixdelindice}
\end{eqnarray}
 By Propositions $ \refe{regularite2} $ and
$ \refe{injectionlinfini}  $, we therefore get, for $ p  \in \{ 2
, \infty \} $,
$$ \left| \left| w(r)^{\frac{n-1}{2} - \frac{n-1}{p}} h^m A (h^2 P - z)^{-1} {\mathcal R}_{N^{\prime}}
(z,h) h^{m^{\prime}} B \right| \right|_{L^2({\mathcal
M},\widetilde{dg})\rightarrow L^{p}({\mathcal M},\widetilde{dg})}
\lesssim \frac{h^{-n/2}}{ |\mbox{Im} \ z|} \left(
\frac{\scal{z}}{|\mbox{Im}\ z |} \right)^{n_2} ,
$$
where the extra power of $ |\mbox{Im}(z)|^{-1} $ comes from the term $ (h^2 P -z)^{-1} $ in (\ref{integrabiliteverticale}). 
Using $ \reff{justificationGreen} $, this estimate clearly proves
that, for $ p \in \{ 2 , \infty \} $,
$$ \left| \left|  w(r)^{\frac{n-1}{2} - \frac{n-1}{p}} h^m A {\mathcal R}_{N^{\prime}}
(P,\varphi,h) h^{m^{\prime}} B \right| \right|_{L^2({\mathcal
M},\widetilde{dg})\rightarrow L^{p}({\mathcal M},\widetilde{dg})}
\lesssim h^{-n(1/2 -1/p)} , $$ if we choose $ M \geq \nu_2 + 1 $ in $
\reff{restedecale} $.
 Then, define $ {\mathcal Q}_{NN^{\prime}}(P,\varphi,h) $ by
\begin{eqnarray}
 \sum_{\iota \in \iota} {\mathcal Q}_{N^{\prime}}^{\iota}
(P,\varphi,h) = \sum_{\iota \in \iota} {\mathcal Q}_N^{\iota}
(P,\varphi,h) + h^{N+1} {\mathcal Q}_{NN^{\prime}}(P,\varphi,h) .
\label{developpementdecale}
\end{eqnarray}
 Using the explicit form of $ {\mathcal
Q}_{NN^{\prime}}(P,\varphi,h) $, namely the fact that its symbol
is a linear combination of terms of the form $ a
(r,\theta,\rho,w(r)\eta) $ with $ a \in S^{-2 \sigma - N} $
(this is due to (\ref{formeexplicitesymbole})), one has
$$ \left| \left|  w(r)^{\frac{n-1}{2} - \frac{n-1}{p}} h^m A {\mathcal Q}_{NN^{\prime}}
(\varphi,h) h^{m^{\prime}} B \right| \right|_{L^2({\mathcal
M},\widetilde{dg})\rightarrow L^{p}({\mathcal M},\widetilde{dg})}
\lesssim h^{-n \left( \frac{1}{2} - \frac{1}{p} \right) } , \qquad
h \in (0,1] ,
$$ which is consequence of Propositions \ref{complocale}, $
\refe{complocale2}  $ and of Lemmas  \ref{appSchur}  and 
\ref{appSchurL2}. Since
\begin{eqnarray}
 {\mathcal R}_{N} (P,\varphi,h) = h^{N^{\prime}-N}{\mathcal
R}_{N^{\prime}} (P,\varphi,h)  + {\mathcal
Q}_{NN^{\prime}}(P,\varphi,h) , \label{doubledeveloppement}
\end{eqnarray}
 by  choosing $ N^{\prime} $ such
that $ N^{\prime} -N - 2 k \geq  - n/ 2 + n / p $ and (\ref{choixdelindice}) holds, we get (\ref{Sobolevnatureltilde})
 for $ p = 2 $ or $ \infty $. The other cases follow by
interpolation. \finpreuve

\section{${ L^p } $ bounds for the resolvent}
\label{preuvetheoremes12}
\setcounter{equation}{0}

Consider a temperate weight $ W $ in the sense of Definition \ref{definitionpoidstempere}.
The main purpose of this section is to prove the following theorem.

\begin{theo} \label{Lptilde}  For
 all $ 1 < p < \infty $, there exists $ \nu_p > 0 $ such that
$$  || W(r) (z-\widetilde{\Delta}_g)^{- 1} W(r)^{-1} ||_{L^p({\mathcal M},\widetilde{dg})}
 \lesssim  \left( \frac{\langle z  \rangle}{|\emph{Im} \ z|} \right)^{\nu_p} , $$
for all $ z \in \Ca \setminus \Ra $.
\end{theo}

Recall that $ \widetilde{\Delta}_g $ is defined by  (\ref{laplacientilde}) and is self-adjoint with respect to $ \widetilde{dg} $ given
by   (\ref{mesuretilde}).

\medskip

Translated in terms of $ \Delta_g $,  Theorem \ref{Lptilde} gives
\begin{coro} For  all $ 1 < p < \infty $, there exists $ \nu_p > 0 $ such that
$$  || W(r) w(r)^{(n-1) (\frac{1}{p}-\frac{1}{2})} (z-\Delta_g)^{- 1} w(r)^{(1-n) ( \frac{1}{p} - \frac{1}{2})} W(r)^{-1} ||_{L^p({\mathcal M},dg)}
 \lesssim \left( \frac{\langle z  \rangle}{|\emph{Im} \ z|} \right)^{\nu_p } , $$
for all $ z \in \Ca \setminus \Ra $.
\end{coro}

Theorem \ref{Lptilde} is a consequence of Proposition \ref{impliqueLprestetilde} showing a stronger result, namely that, in local charts, $ (z-\widetilde{\Delta}_g)^{-1} $
is a pseudo-differential operator with symbol in a class that guarantees the $ L^p $ boundedness on $ L^p ({\mathcal M},\widetilde{dg}) $. Using  Proposition \ref{impliqueLprestetilde}, we also obtain the following result.

\begin{theo} \label{Lpsanstilde} If $ w $ is itself a temperate weight, then for all temperate weight $ W $
and all $
 1 < p < \infty $, there exists $ \nu_p > 0 $ such that
 $$  ||W(r) (z-\Delta_g)^{- 1} W(r)^{-1} ||_{L^p({\mathcal M}, dg)}
 \lesssim  \left( \frac{\langle z  \rangle}{|\emph{Im} \ z|} \right)^{\nu_p}, $$
for all $ z \in \Ca \setminus \Ra $.
\end{theo}
This holds in particular if $ W \equiv 1 $.

\subsection{Reduction}
In this subsection, we explain how to reduce  Theorem \ref{Lptilde} to Proposition 
\ref{impliqueLprestetilde} below. This  reduction rests on  classical results on pseudo-differential
operators, namely the Calder\`on-Zygmund Theorem \ref{MiHo} 
and the Beals Theorem \ref{Beals}.

Recall first the definitions of the usual classes of symbols $ S^0
$ and $ S^0_0 $:
\begin{eqnarray}
a \in S^0 (\Ra^d \times \Ra^d) \Leftrightarrow | \partial_x^{\alpha} \partial_{\xi}^{\beta} a (x,\xi) | \lesssim
\scal{\xi}^{- | \beta |} , \\
a \in S^0_0 (\Ra^d \times \Ra^d) \Leftrightarrow | \partial_x^{\alpha} \partial_{\xi}^{\beta} a (x,\xi) | \lesssim
1 .
\end{eqnarray}
The following theorem is due to Calder\`on-Zygmund.
\begin{theo} \label{MiHo} Let $
d \geq 1 $ and $ a  \in S^0 (\Ra^{d} \times \Ra^{d}) $. Then, for
all $ 1 < p < \infty $,
$$ || a(x,D) v ||_{L^p (\Ra^{d})} \leq C_p || v ||_{L^p (\Ra^{{d}})}, \qquad
v \in C_0^{\infty}(\Ra^{d}) , $$ where the constant $C_p$ depends
on a finite number of semi-norms of $a$ in $ S^0 $.
\end{theo}
For a proof, see for instance \cite{Tayl1}.

We next introduce the class $ S^{-2,0}_{0,1}( \Ra^{n+1} \times
\Ra^n) $ of functions $ b (x_1,x_1^{\prime},y,\rho,\eta) $
satisfying
\begin{eqnarray}
 \left| \partial_{x_1}^j \partial_{x_1^{\prime}}^{j^{\prime}} \partial_{y}^{\alpha} \partial_{\rho}^k
\partial_{\eta}^{\beta}
 b (x_1,x_1^{\prime},y,\rho,\eta) \right| \leq C_{j \alpha k  \beta} \scal{\rho}^{-2} \scal{\eta}^{-|\beta|} ,
  \label{semi-normbis}
\end{eqnarray}
for $ x_1 , x_1^{ \prime} \in \Ra $, $ y \in \Ra^{n-1}$ , and $
(\rho , \eta) \in \Ra \times \Ra^{n-1} $. In particular, for fixed
$ x_1,x_1^{\prime}, \rho $, theses functions belong to $ S^0
(\Ra^{n-1}_{y} \times \Ra^{n-1}_{\eta}) $. 
 Consider the
pseudo-differential operator $B$ defined on $ \Ra^n $ by the  Schwartz kernel 
\begin{eqnarray}
K_B (x_1,y,x_1^{\prime} ,y^{\prime} ) =  (2\pi)^{-n} \int
 e^{ i (y-y^{\prime}) \cdot \eta} \hat{b}
(x_1,x_1^{\prime},y,x_1^{\prime}-x_1,\eta) d\eta \label{lastline}
\end{eqnarray}
where $ \hat{b} $ is the Fourier transform of $b$ with respect to
$ \rho $.  This kernel is continuous with respect to $ x_1 ,
x_1^{\prime} $ (with values in $ {\mathcal
S}^{\prime}(\Ra^{n-1}\times\Ra^{n-1}) $). Integrating by parts
with $ (x_1-x_1^{\prime})^{-1} \partial_{\rho} $ in the integral
defining $\hat{b}$, one sees that, for all $ N $ and all $
\alpha,\beta $,
\begin{eqnarray}
 |\partial_{\theta}^{\alpha} \partial_{\eta}^{\beta} \hat{b}
(x_1,y,x_1^{\prime}-x_1,\eta) | \leq C_{N \alpha \beta}
 \scal{x_1-x_1^{\prime}}^{-N} \scal{\eta}^{-|\beta|} . \label{fastpolynomial}
% \qquad (r , r^{\prime}) \in \Ra \times \Ra, \ (\theta , \eta) \in \Ra^{n-1} \times \Ra^{n-1} .
\end{eqnarray}
Thus, for all $ 1 < p < \infty
$ and $ N > 0 $, Theorem \ref{MiHo} yields the existence of $
C_{Np} $ such that
\begin{eqnarray}
 || (B v)(x_1,.)
||_{L^p(\Ra^{n-1})} \leq
 C_{N p} \int \scal{x_1-x_1^{\prime}}^{-N} || v (x_1^{\prime},.) ||_{L^p(\Ra^{n-1})} d x_1^{\prime}
 , \label{reSch}
\end{eqnarray}
 for all $ v \in C_0^{\infty} (\Ra_{x_1^{\prime}}
\times \Ra^{n-1}_{y^{\prime}}) $. Denoting by  $ p^{\prime}   $ the
conjugate exponent to $ p $, H\"older's
inequality yields
$$ || (B v)(x_1,.) ||_{L^p(\Ra^{n-1})}^p  \lesssim \left( \int \scal{x_1-x_1^{\prime}}^{-N}  d
x_1^{\prime} \right)^{\frac{p}{p^{\prime}}} \left( \int
\scal{x_1-x_1^{\prime}}^{-N} || v (x_1^{\prime},.) ||_{L^p(\Ra^{n-1})}^p
d x_1^{\prime} \right) $$ and thus, if $ N > 1 $, we conclude that
\begin{eqnarray}
|| B v ||_{L^p(\Ra^{n})}^p \lesssim  \int \! \! \int
\scal{x_1-x_1^{\prime}}^{-N} || v (x_1^{\prime},.) ||_{L^p(\Ra^{n-1})}^p
d x_1^{\prime} dx_1 \lesssim || v ||_{L^p (\Ra^n)}^p, \qquad v \in
C_0^{\infty}(\Ra^n) .
\end{eqnarray}
More generally, if $ W $ is a temperate weight, estimates of the form (\ref{fastpolynomial}) still hold if we replace $ \hat{b}
(x_1,y,x_1^{\prime}-x_1,\eta)  $ by $ W (x_1) \hat{b}
(x_1,y,x_1^{\prime}-x_1,\eta) W(x_1^{\prime})^{-1}  $. All this gives the following result.

\begin{prop} \label{MiHoSchu} If $b \in S^{-2,0}_{0,1} (\Ra^{n+1} \times
  \Ra^n)$ and $ B $ is defined by the kernel (\ref{lastline}), then for all temperate weight $W$, $W(x_1) B W(x_1)^{-1} $
  is bounded on $ L^p ( \Ra^{n}) $
for all $ 1 < p < \infty $ and its norm
% $ || W(x_1)^{-1} B W(x_1) ||_{L^p \rightarrowL^p } $
 depends on a finite number of constants $ C_{j\alpha k
\beta} $ in (\ref{semi-normbis}).
\end{prop}

 We shall essentially prove Theorem \ref{Lptilde} by showing that the pull-backs on
 $ \Ra^n $ of $ (z-\widetilde{\Delta}_g)^{-1} $ by local charts are pseudo-differential operators with
 symbols in $ S^{-2,0}_{0,1} (\Ra^{n+1} \times
  \Ra^n) $. The main tool to characterize these pull-packs as
  pseudo-differential operators on $ \Ra^n $ is the
  Beals criterion which we recall in Theorem  \ref{Beals} below.
Fix first some notation.  If $ A $
and $ L $ are operators on suitable spaces,  we set
$$ ad_L \cdot A =  L A - A L .  $$
In our case, $ L $ will typically belong to
% the set of operators
$$ {\mathcal L}_{\Ra^n} = \{ x_1 , \ldots , x_n , \partial_{x_1}, \ldots , \partial_{x_n} \} . $$
\begin{theo}[Beals] \label{Beals} Let $ A : {\mathcal S}(\Ra^n) \rightarrow {\mathcal S}^{\prime}(\Ra^n)
$ be a continuous linear map. If $ A $ is bounded on $ L^2 (\Ra^n)
$ and, more generally for all $ N $ and all $ L_1 , \ldots , L_N
\in {\mathcal L}_{\Ra^n}$, if the operator $ ad_{L_1} \ldots
ad_{L_N} \cdot A  $ is bounded on $ L^2 (\Ra^n) $, then
 there exists $ a \in S^0_0 $  such that
$$ A = a^{W}(x,D) , $$
and each semi-norm of $ a $ in $ S^0_0 $ is controlled by
 a finite number of $ || ad_{L_1} \ldots ad_{L_N} \cdot A  ||_{L^2 \rightarrow L^2} $.
\end{theo}

Here $ a^W (x,D) $ is the Weyl quantization of $ a$ namely the
operator whose kernel is
$$ (2 \pi)^{-n} \int e^{i (x-x^{\prime}) \cdot \xi } a \big( (x+x^{\prime})/2 , \xi \big) d \xi . $$

Theorem \ref{Beals} is for instance proved in \cite{Bea1,Bony,DiSj1}.

The characterization of operators with symbols in $ S^{-2,0}_{0,1}
(\Ra^{n+1} \times \Ra^n) $ is easily deduced from this theorem as
follows. Recall first the formula
\begin{eqnarray}
 (\partial_x^{\alpha} \partial_{\xi}^{\beta} a )^W (x,D) = i^{-|\beta|}
 ad_{\partial_x}^{\alpha} ad_{x}^{\beta} \cdot a^{W}(x,D) ,
\label{deriveesli}
\end{eqnarray}
 where  $ ad_{x}^{\alpha} = ad_{x_1}^{\alpha_1}
\ldots ad_{x_n}^{\alpha_n} $ and $ ad_{\partial_x}^{\beta} =
ad_{\partial_{x_1}}^{\beta_1} \ldots ad_{\partial_{x_n}}^{\beta_n}
$ (note that $ ad_{L_1} ad_{L_2} = ad_{L_2} ad_{L_1} $ for all $
L_1 , L_2 \in {\mathcal L}_{\Ra^n} $). On the other hand, we also
have
\begin{eqnarray}
(\xi_j a)^W(x,D) =  D_j a^W(x,D) - \frac{1}{2 i}
(\partial_{x_j} a)^W (x,D) . \label{deriveeslu}
\end{eqnarray}
\begin{prop} \label{Weylclas}
 Let $ A : {\mathcal S}(\Ra^n) \rightarrow {\mathcal S}^{\prime}(\Ra^n)  $ be linear and continuous. Assume that, for all  $ \alpha, \beta \in \Na^{n}
 $ and all $ \gamma \in \Na^n $ such that
 $$ \gamma_1 \leq 2 , \qquad \gamma_2 + \cdots + \gamma_n \leq \beta_2 + \cdots + \beta_n , $$
 the operator  
\begin{eqnarray}
 A_{ \alpha \beta}^{\gamma} := D_{x}^{\gamma} 
\left(  ad_{\partial_x}^{\alpha}
ad_{x}^{ \beta} \cdot A \right) \label{nomi}
\end{eqnarray}
is bounded on $ L^2 (\Ra^n) $. Then $ A $ is a pseudo-differential operator with symbol  $ a \in S^{-2,0}_{0,1}(\Ra^{n+1} \times \Ra^n)  $ (ie has a kernel of the form (\ref{lastline})). Each semi-norm of $a$ in $
S^{-2,0}_{0,1}(\Ra^{n+1}\times\Ra^n) $ depends on a finite number
of operator norms  $|| A_{ \alpha \beta}^{\gamma} ||_{ L^2
 \rightarrow L^2  } $.
\end{prop}

\noindent {\it Proof.} Set $ B = (1+D_{x_1}^2) A $. By Theorem \ref{Beals}, we can write $ B = b^W (x,D) $ for some $ b \in S^0_0 $. 
Define then $ B^{\gamma}_{\alpha \beta} $ similarly to (\ref{nomi}) with $ B $ instead of $A$ and with $ \gamma = (0,\gamma_2, \ldots , \gamma_n) $.
 By (\ref{deriveesli}) and
 (\ref{deriveeslu}), $ B^{\gamma}_{\alpha \beta} $ is the sum of 
$$ i^{-|\beta|}\left(  \xi^{\gamma}   \partial_{x}^{\alpha} 
\partial_{\xi}^{\beta} b \right)^W (x,D)  , $$
and of a linear combination of operators of the form
$$  \left(  \xi^{\gamma^{\prime}}   \partial_{x}^{\alpha^{\prime}} 
\partial_{\xi}^{\beta} b \right)^W (x,D) , \qquad \gamma^{\prime} < \gamma , \ \ \alpha^{\prime} \leq \alpha + \gamma . $$
On the other hand, by Theorem \ref{Beals} again, $ B^{\gamma}_{\alpha \beta} $ is of the  form $ (b_{\alpha \beta}^{\gamma})^W (x,D) $ for some 
$ b_{\alpha \beta}^{\gamma} \in S_0^0 $. Thus
$$ b_{\alpha \beta}^{\gamma} (x,\xi) = i^{-|\beta|}  \xi^{\gamma}   \partial_{x}^{\alpha} 
\partial_{\xi}^{\beta} b (x,\xi) + \sum_{\gamma^{\prime} < \gamma, \atop \alpha^{\prime} \leq \alpha + \gamma } c_{\gamma^{\prime} \alpha^{\prime}} \xi^{\gamma^{\prime}}  \partial_{x}^{\alpha^{\prime}} 
\partial_{\xi}^{\beta} b (x,\xi) .  $$
By induction on $ \beta $, we deduce that
\begin{eqnarray}
 |  \partial_{x}^{\alpha}
\partial_{\xi}^{\beta} b (x,\xi) | \lesssim  (1+|\xi_2| + \cdots + |\xi_n|)^{-\beta_2 - \cdots - \beta_n} . \label{symbolezero}
\end{eqnarray}
Using then the standard fact that  any $ c^W (y,D_y) $, with $c \in S^0(\Ra^{n-1} \times
\Ra^{n-1})$, can be written $
c_1(y,D_y) $ for some $ c_1 \in S^0
(\Ra^{n-1} \times \Ra^{n-1}) $ depending continuously on $ c $, we can write 
  $ b^W (x,D_x) = b_1 (x,D_x) $ for some symbol $b_1$ satisfying the estimates (\ref{symbolezero}) and depending continuously on $b$.
Therefore $ A = (1+D_{x_1}^2)^{-1} b_1 (x,D_x) $ and its symbol $ \scal{\xi_1}^{-2} b_1 (x,\xi) $ clearly belongs to $ S^{-2,0}_{0,1} $.
\finpreuve

\bigskip

Let us now choose, for each $ \iota \in I $, three functions $ f_{\iota}^{(1)}, f_{\iota}^{(2)},f^{(3)}_{\iota} \in C^{\infty}({\mathcal M}) $ such that, if we set also
$$ f_{\iota}^{(0)} = f_{\iota} $$
$ f_{\iota} $ being the $ \iota $-th element of the partition of unit (\ref{partition}), we have
\begin{eqnarray}
 f_{\iota}^{(j+1)} \equiv 1 \ \ \mbox{near} \ \mbox{supp}(f_{\iota}^{(j)}) , \qquad j = 0, 1, 2 ,
 \label{supportsurcutoff}
\end{eqnarray}
and
\begin{eqnarray}
\mbox{supp} ( f_{\iota}^{(j)} ) \subset {\mathcal U}_{\iota} , \qquad j = 1,2,3 . \label{supportsurutoff2} 
\end{eqnarray}
If $ \iota \in I_{\rm comp} $ we may  assume that $ f_{\iota}^{(j)} \in C_0^{\infty}({\mathcal U}_{\iota}) $ and if $ \iota \in I_{\infty} $
we may assume that 
$$ \Psi_{\iota *} f_{\iota}^{(j)} (r,\theta) = \varrho^{(j)}(r) \kappa_{\iota}^{(j)} (\theta) ,$$ 
with $ \varrho^{(j)} $ and $ \kappa_{\iota}^{(j)} $ supported in small neighborhoods of $ \varrho $ and $ \kappa_{\iota} $ respectively (see (\ref{notationcutoff2})), $ \kappa_{\iota}^{(j)} $ being compactly supported and $ \varrho^{(j)}(r) = 1 $ for $r$ large. Therefore, in all cases,
$$ f_{\iota}^{(j)} \in \mbox{Diff}_w^0 ({\mathcal M}) . $$

By   (\ref{partition}) we can write 
\begin{eqnarray}
(z-\widetilde{\Delta}_g)^{-1} = \sum_{\iota \in I} f_{\iota}^{(0)} (z-\widetilde{\Delta}_g)^{-1} f_{\iota}^{(2)} 
+
\sum_{\iota^{\prime} \in I} \sum_{\iota \in I} f_{\iota}^{(0)} (P-z)^{-1} ( 1-f_{\iota}^{(2)}) f_{\iota^{\prime}}^{(0)} . \nonumber
\end{eqnarray}
The first sum corresponds to `diagonal terms' and the second double one to 'off diagonal terms' since $ f_{\iota}^{(0)} $ and $ ( 1-f_{\iota}^{(2)}) f_{\iota^{\prime}}^{(0)} $ have disjoint supports.

By Proposition \ref{MiHoSchu}, 
Theorem  \ref{Lptilde} would be  a direct consequence of the following proposition.

\begin{prop} \label{impliqueLprestetilde}  For all $ \iota , \iota^{\prime} \in I
$, the operators
\begin{eqnarray}
 R_{\iota} (z) \equiv \Psi_{\iota *} f_{\iota}^{(0)} (z-\widetilde{\Delta}_g)^{-1} f^{(2)}_{\iota} \Psi_{\iota}^* ,
 \qquad z \notin \Ra , \nonumber 
\end{eqnarray}
and
\begin{eqnarray}
R_{\iota \iota^{\prime}} (z) = \Psi_{\iota *}
f_{\iota}^{(0)} (z-\widetilde{\Delta}_g)^{-1} (1 - f^{(2)}_{\iota}) f_{\iota^{\prime}}^{(0)}
\Psi_{\iota^{\prime}}^* , \qquad z \notin \Ra , \nonumber 
\end{eqnarray}
have kernels of the form (\ref{lastline}) with symbols 
whose semi-norms in $ S^{-2,0}_{0,1}(\Ra^{n+1}\times \Ra^n) $ are
bounded by $  ( \scal{z}/| \emph{Im} \
z| )^{\nu} $, for some $\nu$ (depending on the semi-norm).
\end{prop}

We shall prove Proposition \ref{impliqueLprestetilde} using Proposition \ref{Weylclas}. To compute the commutators with elements of $ {\mathcal L}_{\Ra^n} $, we start with a few remarks. For  $ k = 1, \ldots , n $, we have
\begin{eqnarray}
x_k \Psi_{\iota *} = \Psi_{\iota *} x_k^{\iota}, \qquad \Psi_{\iota^{\prime}}^* x_k = x_k^{\iota^{\prime}} \Psi^*_{\iota^{\prime}} ,
\label{coordonnees}
\end{eqnarray}
denoting by $ ( x_1^{\iota} , \ldots , x^{\iota}_n ) $  the coordinates in the $ \iota $-th chart. Similarly
\begin{eqnarray}
\partial_{ x_k } \Psi_{\iota *} = \Psi_{\iota *} \partial_{ x_k^{\iota} }, \qquad \Psi_{\iota^{\prime}}^* \partial_{ x_k } = 
\partial_{ x_k^{\iota^{\prime}} } \Psi^*_{\iota^{\prime}} . \label{vecteur}
\end{eqnarray}
Of course, both (\ref{coordonnees}) and (\ref{vecteur}) hold only in coordinate patches. If $ \iota $ and $ \iota^{\prime} $ belong to $ I_{\infty} $, (\ref{coordonnees}) reads, for
$ k =2 , \ldots , n $,
$$ x_k \Psi_{\iota *} = \Psi_{\iota *} \theta_{k-1}^{\iota}, \qquad \Psi_{\iota^{\prime}}^* x_k = \theta_{k-1}^{\iota^{\prime}} \Psi^*_{\iota^{\prime}} , $$
and for $ k = 1 $,
\begin{eqnarray}
 x_1 \Psi_{\iota *} = \Psi_{\iota *} r, \qquad \Psi_{\iota^{\prime}}^* x_1 = r \Psi^*_{\iota^{\prime}} , \label{pullbackradial}
 \end{eqnarray}
where one should note that $ r $ is globally defined on $ {\mathcal M} $.
We don't write the analogous formulas corresponding to (\ref{vecteur})  for $ \iota, \iota^{\prime} \in I_{\infty} $ but we recall that $ \partial_r $ is only defined where $r$ is a coordinate, namely for $ r > R $.

We then note that $ R_{\iota}(z) $ in Proposition \ref{impliqueLprestetilde} reads
\begin{eqnarray}
 R_{\iota}(z) = \Psi_{\iota *} F_{\iota}^{(0)} R^{\mathcal M} (z) F_{\iota}^{(2)} \Psi^{*}_{\iota} , 
 \label{LPdia}
 \end{eqnarray} 
 with
 \begin{eqnarray}
 \mbox{supp}(F_{\iota}^{(0)}) \subset \mbox{supp}(f_{\iota}^{(0)}), \qquad \mbox{supp}(F_{\iota}^{(2)}) \subset \mbox{supp}(f_{\iota}^{(2)}) , 
 \label{cutoffdiagabstrait} 
 \end{eqnarray}
 and
 \begin{eqnarray}
 R^{\mathcal M} (z) = w(r)^{k_1} A_1 (z-\widetilde{\Delta}_g)^{-1} w(r)^{k_2} A_2 (z-\widetilde{\Delta}_g)^{-1} \cdots w (r)^{k_N} A_N
  (z-\widetilde{\Delta}_g)^{-1}  \label{formeintrinseque}
\end{eqnarray}
 where  $ F_{\iota}^{(0)}, F_{\iota}^{(2)} \in \mbox{Diff}_{w}^{0}({\mathcal M})$,  $ N \geq 1 $, $ k_1 , \ldots , k_N \geq 0 $ and 
$$ A_j \in \mbox{Diff}_{w}^{m_j}({\mathcal M}), \qquad 0 \leq m_j \leq 2 , $$
with in particular $ m_1 = 0 $.  Similarly, 
\begin{eqnarray} 
R_{\iota \iota^{\prime}}(z) = \Psi_{\iota *} F_{\iota}^{(0)} R^{\mathcal M} (z) F_{\iota \iota^{\prime}}^{(2)} \Psi^{*}_{\iota^{\prime}}, \label{LPhd} 
\end{eqnarray}
with $  F_{\iota \iota^{\prime}}^{(2)} \in \mbox{Diff}_{w}^{0}({\mathcal M}) $ such that
\begin{eqnarray}
  \mbox{supp}(F_{\iota \iota^{\prime}}^{(2)}) \subset \mbox{supp} \left( (1-f_{\iota}^{(2)}) f_{\iota^{\prime}}^{(0)} \right) . \label{supporthorsdiagonale} 
\end{eqnarray} 
Of course, in (\ref{LPdia}) and (\ref{LPhd}),  $ R^{\mathcal M}(z) = (z-\widetilde{\Delta}_g)^{-1} $ that is (\ref{formeintrinseque}) with $ N = 1 $, $ k_1 = 0 $ and $A_1 = 1 $. 
 By (\ref{supportsurcutoff}), (\ref{vecteur}) and  (\ref{cutoffdiagabstrait}), we have
\begin{eqnarray}
 \partial_{k} \left( \Psi_{\iota *} F_{\iota}^{(0)} R^{\mathcal M} (z) F_{\iota}^{(2)} \Psi^{*}_{\iota} \right)   = \Psi_{\iota *} \left[ L_{\iota,k} , \ F_{\iota}^{(0)} R^{\mathcal M} (z) F_{\iota}^{(2)}   \right] \Psi_{\iota}^* + \left( \Psi_{\iota *} F_{\iota}^{(0)} R^{\mathcal M} (z) F_{\iota}^{(2)} \Psi^{*}_{\iota}  \right) \partial_k, \nonumber
\end{eqnarray}
with
 \begin{eqnarray}
  L_{\iota,k} = f_{\iota}^{(3)} \Psi_{\iota}^* \partial_{k} \Psi_{\iota *} . \label{champdiagonal}
\end{eqnarray}
In particular,
\begin{eqnarray}
\left[ \partial_{k} ,   \Psi_{\iota *} F_{\iota}^{(0)} R^{\mathcal M} (z) F_{\iota}^{(2)} \Psi^{*}_{\iota}   \right] = \Psi_{\iota *} \left[ L_{\iota,k} ,  F_{\iota}^{(0)} R^{\mathcal M} (z) F_{\iota}^{(2)}  \right] \Psi_{\iota}^* . \label{commutateurdiagonal}
\end{eqnarray}

For operators like the right hand side  of (\ref{LPhd}), 
we use (\ref{supporthorsdiagonale}), that $ f_{\iota}^{(1)} \equiv 1 $ near $ \mbox{supp}(f_{\iota}^{(0)}) $ and that
\begin{eqnarray*}
(1 - f_{\iota}^{(2)}) & \equiv & 0 \qquad \mbox{near} \ \ \mbox{supp}(f_{\iota}^{(1)}) , \\
 (1 - f_{\iota}^{(1)}) & \equiv & \begin{cases} 1 & \mbox{near} \ \ \mbox{supp}(1- f_{\iota}^{(2)}) \\ 
0 & \mbox{near} \ \ \mbox{supp}(f_{\iota}^{(0)}) 
 \end{cases} ,  
\end{eqnarray*}
which follow from (\ref{supportsurcutoff}), to obtain 
\begin{eqnarray}
 \partial_{k} \left( \Psi_{\iota *} F_{\iota}^{(0)} R^{\mathcal M} (z) F_{\iota \iota^{\prime}}^{(2)} \Psi^{*}_{\iota^{\prime}} \right) & = & \Psi_{\iota *} \left[ L_{\iota \rightarrow \iota^{\prime},k} ,  F_{\iota}^{(0)} R^{\mathcal M} (z) F_{\iota \iota^{\prime}}^{(2)}  \right] \Psi_{\iota^{\prime}}^*   , \nonumber \\
  \left( \Psi_{\iota *} F_{\iota}^{(0)} R^{\mathcal M} (z) F_{\iota \iota^{\prime}}^{(2)} \Psi^{*}_{\iota^{\prime}} \right) \partial_k & = & \Psi_{\iota *} \left[   F_{\iota}^{(0)} R^{\mathcal M} (z) F_{\iota \iota^{\prime}}^{(2)}  , L_{\iota \leftarrow \iota^{\prime},k} \right] \Psi_{\iota^{\prime}}^* , \qquad  \nonumber
\end{eqnarray}
with
\begin{eqnarray}
L_{\iota \rightarrow \iota^{\prime},k} = f_{\iota}^{(1)} \Psi_{\iota}^* \partial_{k} \Psi_{\iota *}, \qquad L_{\iota \leftarrow \iota^{\prime},k} = (1-f_{\iota}^{(1)}) f_{\iota^{\prime}}^{(1)} \Psi_{\iota^{\prime} }^* \partial_{k} \Psi_{\iota^{\prime} *} . \label{champhorsdiagonale}
\end{eqnarray}
The main consequence is that
$$
\left[ \partial_{k} , \Psi_{\iota *} F_{\iota}^{(0)} R^{\mathcal M} (z) F_{\iota \iota^{\prime}}^{(2)} \Psi^{*}_{\iota^{\prime}} \right] = \Psi_{\iota *} \left(  \left[ L_{\iota \rightarrow \iota^{\prime},k} ,  F_{\iota}^{(0)} R^{\mathcal M} (z) F_{\iota \iota^{\prime}}^{(2)}  \right]  - 
 \left[  F_{\iota}^{(0)} R^{\mathcal M} (z) F_{\iota \iota^{\prime}}^{(2)}  , L_{\iota \leftarrow \iota^{\prime},k} \right] \right) \Psi_{\iota^{\prime}}^* . 
$$
With the latter formula, (\ref{commutateurdiagonal}) and the resolvent identity, namely 
\begin{eqnarray}
ad_L \cdot (z-\widetilde{\Delta}_g)^{-1} = - (z-\widetilde{\Delta}_g)^{-1}[L,\widetilde{\Delta}_g](z-\widetilde{\Delta}_g)^{-1}, \label{identitedelaresolvente}
\end{eqnarray}
we are equipped to prove the  following result.
\begin{prop} \label{commutateursplusfaciles} For all $ \alpha \in \Na^n $ and all $\iota \in I$ (resp. all $ \iota , \iota^{\prime} \in I $), the operator
$$ ad_{\partial_x}^{\alpha} R_{\iota}(z) \qquad (\mbox{resp.} \ \   ad_{\partial_x}^{\alpha} R_{\iota \iota^{\prime}}(z)), $$
is a linear combination (with coefficients independent of $z$) of operators of the form
$$ F_{\iota}^{(0)} R^{\mathcal M} (z) F_{\iota}^{(2)}  \qquad (\mbox{resp.} \ \ F_{\iota}^{(0)} R^{\mathcal M} (z) F_{\iota \iota^{\prime}}^{(2)} ) $$
with $ N \leq |\alpha| + 1 $ and 
$$ A_1 \in \emph{Diff}_w^0 ({\mathcal M}), \ \  A_2 , \ldots , A_N \in \emph{Diff}_w^2 ({\mathcal M}), \qquad k_1 = k_2 = \cdots = k_N = 0. $$
\end{prop}

\noindent {\it Proof.} It follows from elementary induction once observed that, if $ L $ is any of the operators in (\ref{champdiagonal}) or (\ref{champhorsdiagonale}), we have 
$$ A \in \mbox{Diff}^m_w ({\mathcal M}) \qquad \Rightarrow \qquad \left[L,A \right] \in \mbox{Diff}^m_w ({\mathcal M}) . $$
Indeed, if $ L $ is compactly support this is trivial. Otherwise, if it is supported in chart a infinity,
this is a consequence of the identities
\begin{eqnarray*}
\left[ \partial_r , \widetilde{\varrho}(r) \widetilde{\kappa}(\theta) \partial_{\theta_k} \right] & = & \left(  \frac{\widetilde{\varrho}^{\prime}(r)}{w(r)} \widetilde{\kappa}(\theta) \right) w(r)\partial_{\theta_k} , \\
 \left[ w(r) \partial_{\theta_{k^{\prime}}} , \widetilde{\varrho}(r) \widetilde{\kappa}(\theta) \partial_{\theta_k} \right] & = & \left( \widetilde{\varrho}(r) \partial_{\theta_{k^{\prime}}} \widetilde{\kappa}(\theta)  \right) w(r) \partial_{\theta_k} ,
 % + \left( \widetilde{\varrho}(r) \partial_{\theta_{k^{\prime}}} \widetilde{\kappa}(\theta) \right) w(r) \partial_{\theta_{k^{\prime}}} 
 \\
 \left[  w(r) \partial_{\theta_{k^{\prime}} } , \widetilde{\varrho}(r) \widetilde{\kappa}(\theta) \partial_{r} \right] & = &  \left( w(r)  \widetilde{\varrho}(r) \partial_{\theta_{k^{\prime}}}  \widetilde{\kappa}(\theta) \right)  \partial_r - \left( \frac{w^{\prime}(r)}{w(r)} \widetilde{\varrho}(r) \widetilde{\kappa}(\theta) \right) w (r) \partial_{\theta_{k^{\prime}}} ,
% \left[ w(r) \partial_{\theta_k} , a(r,\theta)  \right] & = & \left( w (r) \partial_{\theta_k} a (r,\theta) \right), 
\end{eqnarray*}
where all the brackets in the right hand sides %the functions $ w^{\prime}(r)/w(r) $, $ w (r) \partial_{\theta_k} a (r,\theta) $ and $ (\varrho^{(j)})^{\prime}(r)/w(r) $ 
are bounded as well as their derivatives, if $ \widetilde{\varrho} $ and $ \widetilde{\kappa} $ are bounded with compactly supported derivatives, also using (\ref{sanscusp}) and (\ref{cinfiniborne}).
\finpreuve
% Notice that all the standard algebraic manipulations are justified by 

\bigskip

To compute $ ad_x^{\beta} ad_{\partial_x}^{\alpha} R_{\iota}(z) $ and $ ad_x^{\beta} ad_{\partial_x}^{\alpha} R_{\iota \iota^{\prime}}(z) $, we need the following lemma.
\begin{lemm} \label{lemmecompactetangulaire} Let $ \widetilde{\rho}  $ be a smooth function on $ \Ra $ with compactly supported derivative and supported in $ r > R $. Let $ \widetilde{\kappa}(\theta) $ be supported in patch of the manifold at infinity.  Then, for any $ A \in \emph{Diff}_w^m ({\mathcal M}) $, we have
$$ \left[A, \widetilde{\varrho}(r) r \right] = A^{\prime} , \qquad  \left[A, \widetilde{\varrho}(r) \widetilde{\kappa}(\theta) \theta_k\right] = w (r) A^{\prime \prime} , $$
for some $ A^{\prime} , A^{\prime \prime} \in \emph{Diff}_w^{m-1} ({\mathcal M}) $. Furthermore, for all $ F \in C_0^{\infty}({\mathcal M}) $ and all $ k \in \Na$, we can write
$$  \left[A, F \right] = w(r)^k A_k , $$
with $ A_k \in \emph{Diff}_w^{m-1} ({\mathcal M}) $.
\end{lemm}

\noindent {\it Proof.} The first two identities follow simply from 
\begin{eqnarray*}
\left[\partial_r , \widetilde{\varrho}(r) r \right] & = & \left( \widetilde{\varrho}^{\prime}(r) r + \widetilde{\varrho}(r)  \right), \\
\left[\partial_r, \widetilde{\varrho}(r) \widetilde{\kappa}(\theta) \theta_k\right] & = & w (r) \left( \frac{\widetilde{\varrho}^{\prime}(r)}{w(r)} \widetilde{\kappa}(\theta) \theta_k  \right) , \\
\left[w(r) \partial_{\theta_{k^{\prime}}}, \widetilde{\varrho}(r) \widetilde{\kappa}(\theta) \theta_k\right] & = & w (r) \widetilde{\varrho}(r) \left( \theta_k \partial_{\theta_{k^{\prime}}} \widetilde{\kappa}(\theta) + \delta_{k k^{\prime}} \widetilde{\kappa}(\theta) \right) ,
\end{eqnarray*}
since all brackets in the right hand sides are smooth and bounded, together with their derivatives. For the third one, we simply observe that  $ \left[A, F \right] $  is a differential operator of order $ m-1 $ with compact support and can thus be written $ w(r)^k (w(r)^{-k} \left[A, F \right] ) $ since $ w  $ doesn't vanish.
\finpreuve

\bigskip

The main sense of this lemma is that commutators of elements of $ \mbox{Diff}_w^{m} ({\mathcal M}) $ with the multiplication operators by coordinates (cut off to be globally defined) are operators in $ \mbox{Diff}_w^{m-1} ({\mathcal M})  $. More precisely, we get a factor $w (r) $ when commuting with angular coordinates or compactly supported functions. Note also that it is crucial for the first commutator that we commute $A$ with a function of $r$ only. Otherwise, we would have to consider for instance terms like
$$ [ w (r) \partial_{\theta_{k^{\prime}}} , \widetilde{\varrho}(r) \widetilde{\kappa}(\theta) r ] = \left( \widetilde{\varrho}(r) \partial_{\theta_{k^{\prime}}} \widetilde{\kappa}(\theta) \right) w(r)r , $$
with $ w(r)r  $ unbounded in general.

To compute the commutators with  the multiplication operators $ x_1 , \ldots , x_n $, we repeat essentially the calculations above Proposition \ref{commutateursplusfaciles} with $ x_k $ instead of $ \partial_k $ except for $ x_1 $ when we work close to infinity. We proceed as follows. If $ \iota \in I_{\rm comp} $, we define
\begin{eqnarray}
 X_{\iota,k} & = &  f_{\iota}^{(3)} \Psi_{\iota}^* x_k \Psi_{\iota *} , \label{X0} \\
X_{\iota \rightarrow \iota^{\prime},k} & = &  f_{\iota}^{(1)} \Psi_{\iota}^* x_k \Psi_{\iota *} , \label{X1}
\end{eqnarray}
for $ 1 \leq k \leq n $. If $ \iota^{\prime } \in I_{\rm comp} $ and $ \iota \in I $, we also set
\begin{eqnarray}
X_{\iota \leftarrow \iota^{\prime},k}  =   (1 - f_{\iota}^{(1)}) f_{\iota^{\prime}}^{(1)}
\Psi_{\iota^{\prime}}^* x_k \Psi_{\iota^{\prime} *} , \label{X2}
\end{eqnarray}
 for $ 1 \leq k \leq n $. In these cases, $ X_{\iota,k} $, $ X_{\iota \rightarrow \iota^{\prime},k} $ and $ X_{\iota \leftarrow \iota^{\prime},k} $ are  smooth functions compactly supported in coordinates patches. If $ k \geq 2 $ and $ \iota, \iota^{\prime} \in I_{\infty} $, we still define $ X_{\iota ,k} $, $ X_{\iota \rightarrow \iota^{\prime},k} $ and $ X_{\iota \leftarrow \iota^{\prime},k} $ by the right hand sides of (\ref{X0}), (\ref{X1})
and (\ref{X2}). Setting finally
\begin{eqnarray}
 X_{\iota,1} & = &  r , \qquad \iota \in I_{\infty} \\
X_{\iota \rightarrow \iota^{\prime},1} & = &   r , \qquad \iota \in I_{\infty}, \ \iota^{\prime} \in I , \\
X_{\iota \leftarrow \iota^{\prime},1} & = & r \qquad \iota \in I , \ \iota^{\prime} \in I_{\infty} ,
\end{eqnarray}
we have defined $ X_{\iota,k} $, $ X_{\iota \rightarrow \iota^{\prime},k} $ and $ X_{\iota \leftarrow \iota^{\prime},k} $ for all $ \iota , \iota^{\prime} \in I $ and all $ 1 \leq k \leq n $. For operators of the form (\ref{formeintrinseque}) and cutoffs satisfying (\ref{cutoffdiagabstrait}), (\ref{coordonnees}) imply that
\begin{eqnarray}
\left[ x_{k} , \left( \Psi_{\iota *} F_{\iota}^{(0)} R^{\mathcal M} (z) F_{\iota }^{(2)} \Psi^{*}_{\iota} \right) \right] & = & \Psi_{\iota *}    F_{\iota}^{(0)} \left[ X_{\iota,k} ,R^{\mathcal M} (z) \right] F_{\iota }^{(2)} \Psi_{\iota}^*   , \nonumber 
\end{eqnarray}
for all $ \iota \in I $ and $ 1 \leq k \leq n $. For off diagonal terms, namely with right cutoffs satisfying (\ref{supporthorsdiagonale}), we have
\begin{eqnarray}
 x_{k}  \left( \Psi_{\iota *} F_{\iota}^{(0)} R^{\mathcal M} (z) F_{\iota \iota^{\prime}}^{(2)} \Psi^{*}_{\iota} \right)  & = & \Psi_{\iota *}    F_{\iota}^{(0)}  X_{\iota \rightarrow \iota^{\prime},k}  R^{\mathcal M} (z)  F_{\iota  \iota^{\prime}}^{(2)} \Psi_{\iota}^*   , \nonumber \\ 
 & = & \Psi_{\iota *}    F_{\iota}^{(0)} \left[ X_{\iota \rightarrow \iota^{\prime},k} ,  R^{\mathcal M} (z) \right] F_{\iota  \iota^{\prime}}^{(2)} \Psi_{\iota^{\prime}}^* \ \ + 
 \nonumber \\
 & & \ \ \qquad \qquad \qquad  \Psi_{\iota *}    F_{\iota}^{(0)} R^{\mathcal M} (z)  X_{\iota \rightarrow \iota^{\prime},k}    F_{\iota  \iota^{\prime}}^{(2)} \Psi_{\iota^{\prime}}^* \nonumber
\end{eqnarray}
where the last term vanishes if $ \iota \in I_{\rm comp} $ or $ k \geq 2 $. In the remaining cases, namely $ k = 1 $ and $ \iota \in I_{\infty} $, we have $ X_{\iota \rightarrow \iota^{\prime},1}  = r $ and
$$  r  F_{\iota  \iota^{\prime}}^{(2)} \Psi_{\iota^{\prime}}^*   = \begin{cases} F_{\iota \iota^{\prime}} \Psi_{\iota^{\prime}}^* \ \ \mbox{with $ F_{\iota \iota^{\prime}} \in C_0^{\infty}({\mathcal M}) $} & \mbox{if} \  \iota^{\prime} \in I_{\rm comp} \ , \\
 F_{\iota  \iota^{\prime}}^{(2)} \Psi_{\iota^{\prime}}^* x_1 & \mbox{if} \  \iota^{\prime} \in I_{\infty}   . \end{cases}  $$
Similarly, we have
\begin{eqnarray}
  \left( \Psi_{\iota *} F_{\iota}^{(0)} R^{\mathcal M} (z) F_{\iota \iota^{\prime}}^{(2)} \Psi^{*}_{\iota} \right) x_k  & = & \Psi_{\iota *}    F_{\iota}^{(0)}    R^{\mathcal M} (z) X_{\iota \leftarrow \iota^{\prime},k}  F_{\iota  \iota^{\prime}}^{(2)} \Psi_{\iota}^*   , \nonumber \\ 
 & = & \Psi_{\iota *}    F_{\iota}^{(0)} \left[   R^{\mathcal M} (z) , X_{\iota \leftarrow \iota^{\prime},k} , \right] F_{\iota  \iota^{\prime}}^{(2)} \Psi_{\iota^{\prime}}^* \ \ + 
 \nonumber \\
 & & \ \ \qquad \qquad \qquad  \Psi_{\iota *}    F_{\iota}^{(0)}   X_{\iota \leftarrow \iota^{\prime},k}  R^{\mathcal M} (z)     F_{\iota  \iota^{\prime}}^{(2)} \Psi_{\iota^{\prime}}^* \nonumber
\end{eqnarray}
where the last term vanishes if $ k \geq 2 $ or $ \iota^{\prime} \in I_{\rm comp} $ and
$$ \Psi_{\iota *}  F_{\iota}^{(0)} X_{\iota \leftarrow \iota^{\prime},1}    =  \begin{cases} F_{\iota \iota^{\prime}} \Psi_{\iota^{\prime}}^* \ \ \mbox{with $ F_{\iota \iota^{\prime}} \in C_0^{\infty}({\mathcal M}) $} & \mbox{if} \  \iota \in I_{\rm comp} \ , \\
x_1 \Psi_{\iota *}  F_{\iota }^{(0)}  & \mbox{if} \  \iota \in I_{\infty}   . \end{cases}  $$

This shows that, unless $ \iota, \iota^{\prime} \in I_{\infty} $ and $ k = 1 $, $  \left[ x_k , \left( \Psi_{\iota *} F_{\iota}^{(0)} R^{\mathcal M} (z) F_{\iota \iota^{\prime}}^{(2)} \Psi^{*}_{\iota} \right) \right] $ is the sum of
$$ \Psi_{\iota *}    F_{\iota}^{(0)} \left(  \left[ X_{\iota \rightarrow \iota^{\prime},k} ,  R^{\mathcal M} (z) \right] - \left[   R^{\mathcal M} (z) , X_{\iota \leftarrow \iota^{\prime},k} , \right] \right) F_{\iota  \iota^{\prime}}^{(2)} \Psi_{\iota^{\prime}}^*   $$
and of terms of the same form as $ \Psi_{\iota *} F_{\iota}^{(0)} R^{\mathcal M} (z) F_{\iota \iota^{\prime}}^{(2)} \Psi^{*}_{\iota} $. If $ \iota, \iota^{\prime} \in I_{\infty} $ and $ k = 1 $, we simply have
$$ \left[ x_{1} , \left( \Psi_{\iota *} F_{\iota}^{(0)} R^{\mathcal M} (z) F_{\iota \iota^{\prime}}^{(2)} \Psi^{*}_{\iota} \right) \right]   =   \Psi_{\iota *}    F_{\iota}^{(0)} \left[ r ,  R^{\mathcal M} (z) \right] F_{\iota  \iota^{\prime}}^{(2)} \Psi_{\iota^{\prime}}^* . $$
Using lemma \ref{lemmecompactetangulaire}, the resolvent identity (\ref{identitedelaresolvente}) and a simple induction, we get the following result.

\begin{prop} \label{algebrenonexplicitee} For all $ \alpha,\beta \in \Na^n $ and all $\iota \in I$ (resp. all $ \iota , \iota^{\prime} \in I $), the operator
$$ ad_x^{\beta} ad_{\partial_x}^{\alpha} R_{\iota}(z) \qquad (\mbox{resp.} \ \  ad_x^{\beta} ad_{\partial_x}^{\alpha} R_{\iota \iota^{\prime}}(z)), $$
is a linear combination (with coefficients independent of $z$) of operators of the form
$$ F_{\iota}^{(0)} R^{\mathcal M} (z) F_{\iota}^{(2)}  \qquad (\mbox{resp.} \ \ F_{\iota}^{(0)} R^{\mathcal M} (z) F_{\iota \iota^{\prime}}^{(2)} ) $$
where $ R^{\mathcal M} (z)  $ is of the form (\ref{formeintrinseque}) with $ N \leq |\alpha| + |\beta| + 1 $,
$$ A_j \in \emph{Diff}_w^{m_j} ({\mathcal M}), \qquad 0 \leq m_j \leq 2 $$
and
$$ k_j = 2 - m_j, \qquad k_2 + \cdots + k_N = \beta_2 + \cdots + \beta_n . $$
\end{prop}

\bigskip

The next proposition is  the final step before being in position to use Proposition \ref{Weylclas}.

\begin{prop} \label{etapefinale} Fix $\iota \in I$ (resp. $\iota, \iota^{\prime} \in I $). For all $ \alpha, \beta \in \Na^n $ and all $ \gamma \in \Na $ satisfying $ \gamma_1 \leq 2 $, $ \gamma_2 + \cdots + \gamma_n \leq \beta_2 + \cdots + \beta_n $, the operator
$$ D_x^{\gamma} ad_x^{\beta} ad_{\partial_x}^{\alpha} R_{\iota}(z) \qquad (\mbox{resp.} \ \ D_x^{\gamma} ad_x^{\beta} ad_{\partial_x}^{\alpha} R_{\iota \iota^{\prime}}(z)), $$
is a linear combination (with coefficients independent of $z$) of operators of the form
$$ F_{\iota}^{(0)} R^{\mathcal M} (z) F_{\iota}^{(2)}  \qquad (\mbox{resp.} \ \ F_{\iota}^{(0)} R^{\mathcal M} (z) F_{\iota \iota^{\prime}}^{(2)} ) $$
(see (\ref{formeintrinseque})) with $ N \leq |\alpha| + |\beta| + |\gamma| + 1 $ and 
$$ A_1 , \ldots , A_N \in \emph{Diff}_w^{2} ({\mathcal M}), k_1 = \cdots = k_N = 0 . $$
In particular, they are bounded on $ L^2 (\Ra^n) $ with norms controlled by powers of $ \scal{z}/|\emph{Im}(z)| $.
\end{prop}

\noindent {\it Proof.} We treat the case of $ R_{\iota}(z) $, the one of $ R_{\iota   \iota^{\prime}}(z) $ being completely similar. We start with a simple case. Consider an operator of the form
$$ B(z) := \Psi_{\iota *} F_{\iota}^{(0)} (z - \widetilde{\Delta}_g)^{-1} w(r) A (z - \widetilde{\Delta}_g)^{-1} F_{\iota}^{(2)} \Psi_{\iota}^* , $$
with $A \in \mbox{Diff}_w^1 ({\mathcal M}) $. Such operators appear in Proposition \ref{algebrenonexplicitee} with $ N = 2 $ if $ \beta_2 + \cdots + \beta_n = 1 $
and $ \alpha = 0 $. Compute then $ \partial_{k} B (z) $, with $k \geq 2$. We get
$$ \Psi_{\iota *} \left( \left[ L_{\iota,k} , F_{\iota}^{(0)} (z - \widetilde{\Delta}_g)^{-1} w(r) \right] A (z - \widetilde{\Delta}_g)^{-1} +  F_{\iota}^{(0)} (z - \widetilde{\Delta}_g)^{-1} w(r)  L_{\iota,k} A (z - \widetilde{\Delta}_g)^{-1} F_{\iota}^{(2)} \right) \Psi_{\iota}^* . $$
The commutator reads
$$ \left[ L_{\iota,k} , F_{\iota}^{(0)} \right] (z - \widetilde{\Delta}_g)^{-1} w(r) +  F_{\iota}^{(0)} (z - \widetilde{\Delta}_g) \left[ \widetilde{\Delta}_g , L_{\iota,k} \right] (z - \widetilde{\Delta}_g)^{-1} w(r)  +  F_{\iota}^{(0)} (z - \widetilde{\Delta}_g)^{-1} \left[L_{\iota,k} ,  w(r) \right] $$
and is bounded on $ L^2 ({\mathcal M},\widetilde{dg}) $ since $ [L_{\iota,k}, \widetilde{\Delta}_g] \in \mbox{Diff}_w^1 ({\mathcal M}) $. The simple and crucial remark is that
$$ w(r)  L_{\iota,k} A \in \mbox{Diff}_w^2 ({\mathcal M}) ,$$
although $ L_{\iota,k} A \notin \mbox{Diff}_w^2 ({\mathcal M}) $ in general. Therefore $ \partial_k B (z) $ is a linear combination of operators of the form (\ref{formeintrinseque}) with $ A_1 $ of order $ 0 $. This then implies that 
$ \partial_1^2 \partial_k B (z) $ is also of this form with $A_1 $ of order 2. Iteration of this argument  give the result  since  Proposition \ref{algebrenonexplicitee} show there are at least $ \gamma_2 + \cdots + \gamma_n $ powers $w(r) $ in the expression of $ ad_x^{\beta} ad_{\partial_x}^{\alpha} R_{\iota}(z) $ to absorb $ \partial_{x_2}^{\gamma_2} \ldots \partial_{x_n}^{\gamma_n} $. \finpreuve

\subsection{Proof of Proposition \ref{impliqueLprestetilde}}
 This is a direct consequence of Proposition \ref{etapefinale} and Proposition \ref{Weylclas}.

\subsection{Proof of Theorem \ref{Lptilde}}
This is a direct consequence of Proposition \ref{MiHoSchu} and Proposition \ref{impliqueLprestetilde} using the equivalence of norms (\ref{equivalenceLpMRn}).
\subsection{Proof of Theorem \ref{Lpsanstilde}}
 The boundedness of $
W(r)(z-\Delta_g)^{-1}W(r)^{-1} $ on $ L^p ({\mathcal M},dg) $ is equivalent
to the one of
\begin{eqnarray}
 W(r) w(r)^{\frac{n-1}{2}-\frac{n-1}{p} } (z-\widetilde{\Delta}_g)^{-1} w(r)^{\frac{n-1}{p}-\frac{n-1}{2}}W(r)^{-1} \nonumber
\end{eqnarray}
  on
$ L^p ({\mathcal M}, \widetilde{dg}) $ so the result follows  from Proposition \ref{MiHoSchu}, with the temperate weight $ W w^{\frac{n-1}{2}-\frac{n-1}{p}} $, and Proposition \ref{impliqueLprestetilde}.
 \finpreuve

\subsection{Proof of Theorem \ref{Lprestetilde}}
We note first that, by writing $ ( z - h^2 \widetilde{\Delta}_g )^{-1} = h^{-2} (zh^{-2} - \widetilde{\Delta}_g)^{-1} $, Theorem \ref{Lptilde} implies that
\begin{eqnarray}
|| W(r) (z-h^{2} \widetilde{\Delta}_g)^{-1} W(r)^{-1} ||_{L^p ({\mathcal M},\widetilde{dg}) \rightarrow L^p ({\mathcal M},\widetilde{dg}) }
\lesssim h^{-2} \frac{\scal{z}^{\nu_p}}{|\mbox{Im}(z)|^{\nu_p}} , \qquad h \in (0,1] , \ \ z \in \Ca \setminus \Ra , \label{estimationsemiclassique}
\end{eqnarray}
by using the inequality $ \scal{h^{-2}z}/{|\mbox{Im}(h^{-2}z)|} \lesssim  \scal{z}/|\mbox{Im}(z)| $.

Assume next that $ \varphi \in S^{-\sigma}(\Ra) $ with $ \sigma > 1 $ so that
we can use (\ref{justificationGreenbis}).
% Then, if $ M \geq \nu  $ in (\ref{dbarre}), we have 
%\begin{eqnarray}
%\int \! \! \int_{\Ra^2} | \bar{\partial} \widetilde{\varphi}_M
%(x+iy) | \times  \left( \frac{ 1+|x|+ |y| }{|y|}
%\right)^{\nu} dx dy < \infty . \label{justificationGreenbis}
%\end{eqnarray}
By  Theorem \ref{parametrixegenerale} and (\ref{estimationsemiclassique}), there exists $ \nu_{p,N} $
such that
$$    \left| \left|W(r) (z-h^2 \widetilde{\Delta}_g)^{-1} {\mathcal R}_{N} (z,h) W(r)^{-1}
  \right| \right|_{L^p ({\mathcal M},\widetilde{dg})\rightarrow
L^p ({\mathcal M},\widetilde{dg}) } \lesssim h^{-2}
\left(
\frac{\scal{z}}{|\mbox{Im} \ z|} \right)^{\nu_p + \nu_{N,p}} , 
$$
for $  h \in (0,1 ] $ and $ z \notin \Ra  $. By choosing $ M \geq \nu = \nu_p
+ \nu_{N,p} $, the above estimate and (\ref{justificationGreenbis}) give the
expected estimate up to a factor $ h^{-2}  $. The latter is eliminated in the
standard way: by pushing the expansion to the order $h^{N+2}$, we write  $ {\mathcal R}_{N} (- \widetilde{\Delta}_g , \varphi ,
h) $ as the sum of 
properly supported pseudo-differential operators bounded on $ W (r)^{-1} L^p
({\mathcal M},\widetilde{dg})  $ and of $ h^2  {\mathcal R}_{N+2} (- \widetilde{\Delta}_g , \varphi ,
h)  $. This implies (\ref{Lppourlerestetilde}).

If now $  \varphi \in S^{-\sigma}(\Ra) $ with $ \sigma > 0  $, we cannot use
(\ref{justificationGreenbis}). We thus  
write $ \varphi (\lambda) = (\lambda + i) \psi (\lambda) $ with $ \psi \in S^{- \sigma -1}(\Ra) $ so that
\begin{eqnarray}
 \varphi (-h^2 \widetilde{\Delta}_g) = ( i - h^2 \widetilde{\Delta}_g) \psi (-h^2 \widetilde{\Delta}_g) . \label{trick}
\end{eqnarray}
We then write again 
$ {\mathcal R}_{N} (- \widetilde{\Delta}_g , \varphi , h) $ as a finite sum of properly supported pseudo-differential operators bounded on $ W (r)^{-1} L^p ({\mathcal M},\widetilde{dg})  $  and 
$$ h^{N+2} \int \! \! \int_{\Ra^2}  \bar{\partial} \widetilde{\psi}_M
(z)  (z-h^2 \widetilde{\Delta}_g )^{-1} (i-h^2 \widetilde{\Delta}_g) {\mathcal R}_{N+2} (z,h) dx dy $$
where $ z = x + i y $. By Theorem \ref{parametrixegenerale}, we have
$$    \left| \left| W(r) (i-h^2 \widetilde{\Delta}) {\mathcal R}_{N+2} (z,h) W(r)^{-1}
  \right| \right|_{L^p ({\mathcal M},\widetilde{dg})\rightarrow
L^p ({\mathcal M},\widetilde{dg}) } \lesssim
\left(
\frac{\scal{z}}{|\mbox{Im} \ z|} \right)^{\nu_{N+2,p}}  ,
$$
and we proceed as above. \finpreuve

\appendix

\section{Non $ L^p \rightarrow L^p $ boundedness on the hyperbolic space} \label{contreexemple}
Using the hyperboloid model of the hyperbolic space, namely
$$ {\mathbb H}^n = \{ x = (x_0, \ldots ,x_n) \in \Ra^{n+1} \ | \ x_0^2 - x_1^2 - \cdots - x_n^2 = 1 , \ x_0 > 0 \} , $$
we have polar coordinates by considering
$$ x (r,\omega) = ( \cosh r , \omega \sinh r ), \qquad r > 0 , \ \omega \in {\mathbb S}^{n-1} . $$
In this parametrization, the distance between $ x = x (r,\omega) $ and $ x^{\prime} = x (r^{\prime},\omega^{\prime}) $ reads
\begin{eqnarray}
 d (x,x^{\prime}) & = & \mbox{arccosh} \left( \cosh r \cosh r^{\prime} - \omega \cdot \omega^{\prime} \sinh r \sinh r^{\prime} \right) 
 \nonumber \\
 & = & \mbox{arccosh} \left\{ \left(1- \frac{|\omega-\omega^{\prime}|^2}{4} \right) \cosh (r - r^{\prime}) + \frac{|\omega-\omega^{\prime}|^2}{4} \cosh ( r + r^{\prime} ) \right\} \label{formeautiliserpourlescalculs}
\end{eqnarray}
and the volume element is
$$ (\sinh r)^{n-1}  dr d \omega , $$
where $ d \omega $ is the usual Riemannian measure on the sphere. Considering $ n = 3 $ for simplicity,  the resolvent
\begin{eqnarray}
 ( - \Delta_{{\mathbb H}^3} - 1 + \epsilon^2 )^{-1} , \qquad \epsilon > 0 , \label{resolventeexplicite}
 \end{eqnarray}
is well defined since, in general, $ - \Delta_{{\mathbb H}^n} \geq (n-1)^2/4 $. Its kernel
with respect to the volume element is then given by
$$ \frac{1}{4 \pi} \frac{e^{- \epsilon d (x,x^{\prime})}}{\sinh d (x,x^{\prime})} . $$ 
 (see for instance \cite[p. 105]{Tayl2}).

\begin{prop} \label{propositionnonbornee} Fix $ p \in (1,\infty) $ with $ p \ne 2 $. If $ 0 < \epsilon < \left| 1 - \frac{2}{p} \right| $, 
then $ ( - \Delta_{{\mathbb H}^3} - 1 + \epsilon^2 )^{-1} $ is not bounded on $ L^p ({\mathbb H}^3) $.
\end{prop}
We shall proceed by contradiction, using the following simple lemma.
\begin{lemm} \label{principedumaximum} Let $ K_1 , K_2 $ be two locally integrable functions on $ ( \Ra_+ \times {\mathbb S}^2 )^2 $ such that
\begin{eqnarray}
 K_2 (r,\omega,r^{\prime},\omega^{\prime}) \geq |K_1 (r,\omega,r^{\prime},\omega^{\prime}) | . \label{noyauxordonnes}
\end{eqnarray} 
Denote by $ A_j $ be the operator with kernel $ K_j $ with respect to $ d r d \omega $ and set $ L^p = L^p (\Ra_+ \times {\mathbb S}^2,drd\omega) $. Then
$$ || A_1 ||_{L^p \rightarrow L^p} \leq || A_2 ||_{L^p \rightarrow L^p} . $$
% \ \mbox{bounded on} \ L^p (\Ra_+ \times {\mathbb S}^2,drd\omega) \qquad \Rightarrow \qquad  A_1 \ \mbox{bounded on} \ L^p (\Ra_+ \times {\mathbb %S}^2,drd\omega) . $$
\end{lemm}

\noindent {\it Proof.} By (\ref{noyauxordonnes}), we have, for all $ u \in C_0^{\infty}( \Ra_+ \times {\mathbb S}^2 ) $, 
$$ |(A_1 u)(r,\omega) | %= \left| \int K_1 (r,\omega,r^{\prime},\omega^{\prime})u(r^{\prime},\omega^{\prime})
%d r^{\prime} d \omega^{\prime} \right| \leq  \int K_2 (r,\omega,r^{\prime},\omega^{\prime})|u(r^{\prime},\omega^{\prime})|
%d r^{\prime} d \omega^{\prime} = 
\leq |(A_2|u|)(r,\omega)| $$ so, taking the $ L^p $ norm, we obtain 
$$ ||A_1 u ||_{L^p} \leq \big| \big| A_2 |u| \big| \big|_{L^p} \leq || A_2 ||_{L^p \rightarrow L^p} \big| \big| |u| \big| \big|_{L^p} = ||A_2 ||_{L^p \rightarrow L^p} ||u||_{L^p} $$
which gives the result. \finpreuve
% the right hand side is bounded by $ C \big| \big| |u| \big| \big|_{L^p} = C ||u||_{L^p} $ if $ A_2 $ is bounded on $ L^p $. \finpreuve

\bigskip

\noindent {\it Proof of Proposition \ref{propositionnonbornee}.} We argue by contradiction and assume that $ ( - \Delta_{{\mathbb H}^3} - 1 + \epsilon^2 )^{-1} $ is bounded on $ L^p ({\mathbb H}^3) $. This is equivalent to the boundedness on $ L^p (\Ra_+ \times {\mathbb S}^2,drd\omega) $ of the operator with kernel
$$ K_2 (r,\omega,r^{\prime},\omega^{\prime}) := (\sinh r)^{\frac{2}{p}} \left( \frac{1}{4 \pi} \frac{e^{- \epsilon d (x,x^{\prime})}}{\sinh d (x,x^{\prime})} (\sinh r^{\prime})^2 \right) (\sinh r^{\prime})^{-\frac{2}{p}}  $$
with respect to $ d r d \omega $. Since $ \cosh (r - r^{\prime}) \leq \cosh (r+r^{\prime}) $ for $ r,r^{\prime} \in \Ra^+ $, (\ref{formeautiliserpourlescalculs})  gives
$$ d (x,x^{\prime}) \leq r + r^{\prime} $$
so, for $ r,r^{\prime} \geq 1 $, we have 
\begin{eqnarray}
 K_2 (r,\omega,r^{\prime},\omega^{\prime}) \gtrsim  (e^r)^{\frac{2}{p}} \left( \frac{e^{- \epsilon (r+r^{\prime})}}{e^{ r + r^{\prime}}} (e^{ r^{\prime}})^2 \right) (e^{r^{\prime}})^{-\frac{2}{p}}  = e^{ (\frac{2}{p}-1 -\epsilon)r} e^{(1- \frac{2}{p}- \epsilon)r^{\prime}} . \label{nontempere}
\end{eqnarray}
Denoting by $ K_1 (r,\omega,r^{\prime},\omega^{\prime}) = K_1 (r,r^{\prime}) $ the right hand side of (\ref{nontempere}) multiplied by the characteristic function of $ [1,+\infty )^2 $, Lemma \ref{principedumaximum} implies that the corresponding operator $ A_1 $ is bounded on $ L^p (\Ra_+ \times {\mathbb S}^2,drd\omega) $. This is clearly not true if $ \frac{2}{p} -1 > \epsilon $, otherwise $ e^{ (\frac{2}{p}-1 -\epsilon)r} $ should belong to $ L^p (\Ra) $. We also obtain a contradiction if $ 1 - \frac{2}{p} > \epsilon $ by considering the adjoint of $A_1 $ \finpreuve

\bigskip

Note that the right hand side of  (\ref{nontempere}) also reads
$$ e^{ (\frac{2}{p}-1)(r-r^{\prime}) -\epsilon ( r + r^{\prime})} , $$
showing that the above reasoning  gives no contradiction  for $ p = 2 $ nor by restricting the kernel close to the diagonal.
% ie to $ |r-r^{\prime}| \lesssim 1 $, the right hand side of (\ref{nontempere}) is of order
%$$  $$

We also recall that  $ (n-1)|\frac{1}{p}-\frac{1}{2}| $ (ie $ |\frac{2}{p}-1| $ if $n=3$) is exactly the width of the strip
around the real axis in which $ \varphi $ has to be holomorphic to ensure the boundedness on $ L^p ({\mathbb H}^n) $ of
$$ \varphi \left( ( - \Delta_{{\mathbb H}^n} - (n-1)^2/4)^{1/2} \right) , $$
as proved in \cite{Tayl0}. The resolvent (\ref{resolventeexplicite}) corresponds to $ \varphi (\lambda) = (\lambda^2 + \epsilon^2)^{-1} $
which is holomorphic for $ |\mbox{Im}(\lambda)| < \epsilon $.

\end{document}